\documentclass[10pt, a4paper]{article}
\usepackage{amsmath,amsthm,amssymb}
\usepackage[plain]{fullpage}

\usepackage{tikz}
\usepackage{color,hyperref}
\usepackage{comment}
\usepackage{authblk}
\usepackage{enumitem}
\usepackage{thm-restate}

\usepackage{tikz, tkz-graph, tkz-berge}
\usetikzlibrary{calc}
\definecolor{g-green}{rgb}{0.235, 0.659, 0.322}
\definecolor{g-blue}{rgb}{0.0, 0.5, 1.0}

\newtheorem{theorem}{Theorem}

\newtheorem{proposition}[theorem]{Proposition}
\newtheorem{lemma}[theorem]{Lemma}
\theoremstyle{definition}

\newtheorem{problem}[theorem]{Problem}
\newtheorem{conjecture}[theorem]{Conjecture}

\newtheorem{claim}{Claim}[theorem]

\newenvironment{proofclaim}[1][]
	{\par\noindent {\it Proof of claim}. }{ \hfill$\lozenge$\par\vspace{11pt}}

\newcommand{\mader}{\mathrm{mader}}

\newcommand{\bid}[1]{\overleftrightarrow{#1}}
\newcommand{\ori}[1]{\overrightarrow{#1}}
\newcommand{\dic}{\vec{\chi}}
\newcommand{\digirth}{\mathrm{digirth}}
\newcommand{\girth}{\mathrm{girth}}
\newcommand{\init}{\mathrm{init}}
\newcommand{\term}{\mathrm{term}}
\newcommand{\pred}{\mathrm{pred}}
\renewcommand{\succ}{\mathrm{succ}}
\renewcommand{\mod}{\mathrm{~mod~}}
\newcommand{\length}{\mathrm{length}}
\newcommand{\ind}[1]{[#1]}

\DeclareMathOperator{\UG}{UG}

\DeclareMathOperator{\dist}{dist}
\DeclareMathOperator{\cc}{cc}

\begin{document}

\title{Subdivisions in dicritical digraphs with large order or digirth\thanks{Research supported by the research grant
    DIGRAPHS ANR-19-CE48-0013 and by the French government, through the EUR DS4H Investments in the Future project managed by the National Research Agency (ANR) with the reference number ANR-17-EURE-0004.}
    }

\author{Lucas Picasarri-Arrieta}
\author{Clément Rambaud}
\affil{Universit\'e C\^ote d'Azur, CNRS, Inria, I3S, Sophia-Antipolis, France}

\date{}
\maketitle
\vspace{-10mm}
\begin{center}
{\small 
\texttt{$\{$lucas.picasarri-arrieta, clement.rambaud$\}$@inria.fr}\\
}
\end{center}

\maketitle

\sloppy

\begin{abstract}
Aboulker et al. proved that a digraph with large enough dichromatic number contains any fixed digraph as a subdivision.
The dichromatic number of a digraph is the smallest order of a partition of its vertex set into acyclic induced subdigraphs.
A digraph is dicritical if the removal of any arc or vertex decreases its dichromatic number.
In this paper we give sufficient conditions on a dicritical digraph of large order or large directed girth to contain a given digraph as a subdivision.
In particular, we prove that 
(i) for every integers $k,\ell$, large enough dicritical digraphs with dichromatic number $k$ contain an orientation of a cycle with at least $\ell$ vertices;
(ii) there are functions $f,g$ such that for every subdivision $F^*$ of a digraph $F$, digraphs with directed girth at least $f(F^*)$ and dichromatic number at least $g(F)$ contain a subdivision of $F^*$, and if $F$ is a tree, then $g(F)=|V(F)|$;
(iii) there is a function $f$ such that for every subdivision $F^*$ of $TT_3$ (the transitive tournament on three vertices), digraphs with directed girth at least $f(F^*)$ and minimum out-degree at least $2$ contain $F^*$ as a subdivision.\\

\noindent{}{\bf Keywords:} Digraphs, dichromatic number, dicritical digraphs, subdivisions, digirth.
\end{abstract}

\section{Introduction}

Since the seminal works of Mader~\cite{maderMA174, maderCOMB5},
a lot of sufficient conditions for a (di)graph to contain a subdivision of a given (di)graph have been proven,
in particular having large chromatic number~\cite{bondyJLMS2, cohenJGT89, kimJGT88}, or having large minimum (out-)degree~\cite{alonJCTB68, bollobasComb16, gishbolinerComb42, Thomassen1983}.
In this paper, we give sufficient conditions on a digraph of large dichromatic number or large out-degree to contain a given digraph as a subdivision.

Let $D$ be a digraph. A {\it $k$-colouring} of $D$ is a function $\phi\colon V(D) \to [k]$. It is a {\it $k$-dicolouring} if no directed cycle $C$ in $D$ is monochromatic for $\phi$. Equivalently, it is a $k$-dicolouring  if every colour class induces an acyclic subdigraph.
The smallest integer $k$ such that $D$ has a $k$-dicolouring is the {\it dichromatic number} of $D$ and is denoted by $\dic(D)$.
For digraphs of large dichromatic number, the most general result is the following.

\begin{theorem}[Aboulker et al.~\cite{aboulkerEJC26}]
    \label{thm:all_digraphs_dic_maderian}
    Let $F$ be a digraph on $n$ vertices, $m$ arcs and $c$ connected components. Every digraph $D$ satisfying $\dic(D) \geq 4^{m - n +c}(n-1) + 1$ contains a subdivision of $F$.
\end{theorem}

For every digraph $F$, we denote by $\mader_{\dic}(F)$ the least integer $c_F$ for which every digraph $D$ with dichromatic number $\dic(D)\geq c_F$ contains a subdivision of $F$. Note that $\mader_{\dic}(F)$ is well-defined by Theorem~\ref{thm:all_digraphs_dic_maderian}.
The result above was generalized in a recent work of Steiner~\cite{steiner2024subdivisions} (see also~\cite{kwon2023variant}) who extended it to subdivisions with modular constraints.

Since every digraph is a subdigraph of $\bid{K_n}$, the complete digraph on $n$ vertices, it is natural to look for the value of $\mader_{\dic}(\bid{K_n})$. The result above implies that $\mader_{\dic}(\bid{K_n}) \leq 4^{n^2-2n+1}$. 
A more precise computation using the tools developed in~\cite{aboulkerEJC26} yields $\mader_{\dic}(\bid{K_n}) \leq 4^{\frac{2}{3}n^2 + 2n - \frac{8}{3}}$, as we show in Section~\ref{section:mader_Kn}.

For every digraph $F$ on $n$ vertices, we have $\mader_{\dic}(F) \geq n$. This is because $\bid{K_{n-1}}$ has dichromatic number $n-1$ and does not contain any subdivision of $F$. 
A digraph $D$ is $k$-dicritical if $\dic(D) = k$ and every proper subdigraph $H$ of $D$ satisfies $\dic(H) < k$. 
For some digraphs $F$, it then appears that the value of $\mader_{\dic}(F)$ does not capture the structure of $F$ but only its order. For instance, $\mader_{\dic}(\ori{C_n}) \geq n$, but every $2$-dicritical 
digraph on at least $n$ vertices actually contains a subdivision of $\ori{C_n}$. In order to have a better understanding of digraphs forced to contain subdivisions of $F$, one may then ask for the minimum $k$ such that there is a finite number of $k$-dicritical digraphs which do not contain any subdivision of $F$. The following question then naturally arises.
    Let $F^*$ be a subdivision of a digraph $F$, is it true that the set of $(\mader_{\dic}(F))$-dicritical digraphs that do not contain any subdivision of $F^*$ is finite?

\begin{figure}[hbtp]
        \begin{center}	
              \begin{tikzpicture}[thick,scale=1, every node/.style={transform shape}]
                \tikzset{vertex/.style = {circle,fill=black,minimum size=4pt, inner sep=0pt}}
                \tikzset{edge/.style = {->,> = latex'}}
                
                \node[vertex] (x) at (0,0) {};
                
                 \foreach \i in {0,...,8}{
                    \node[vertex] (u\i) at (40*\i:2) {};
                    \draw[edge, bend left=12] (u\i) to (x) {};
                    \draw[edge, bend left=12] (x) to (u\i) {};
          	}
                \foreach \i in {0,...,8}{
                    \pgfmathtruncatemacro{\j}{Mod(\i + 1,9)}
                    \draw[edge, bend right=13] (u\i) to (u\j) {};
          	}
              \end{tikzpicture}
          \caption{The digraph $D_{10}$.}
          \label{fig:D_n}
        \end{center}
    \end{figure}
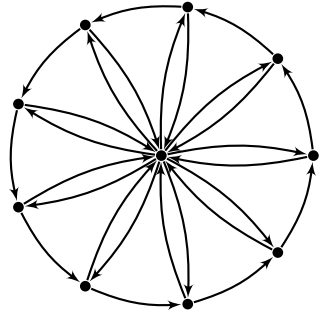

Unfortunately, the answer to this question is negative. 
To see that, consider for every positive integers $k$ and $\ell$ the digraph $C(k,\ell)$, which is the union of two internally disjoint directed paths from a vertex $x$ to a vertex $y$ of lengths respectively $k$ and $\ell$.
Observe that $\mader_{\dic}(C(1,2)) = 3$ because $\bid{K_2}$ does not contain any subdivision of $C(1,2)$ and every $3$-dicritical digraph is $2$-arc-strong. 
However, for every integer $n$ with $n \geq 3$, the digraph obtained from a directed cycle on $n-1$ vertices $\ori{C_{n-1}}$ by adding a new vertex $x$ and all possible digons between $x$ and $V(\ori{C_{n-1}})$ (see Figure~\ref{fig:D_n} for an illustration) is $3$-dicritical (see Lemma~\ref{lemma:universalvertex}) but does not contain any subdivision of $C(3,3)$.

In fact, as we show in Section~\ref{section:large_dicritical}, the answer to the aforementioned question is negative not only for $C(1,2)$ but for every digraph $F$ on at least three vertices with at least one arc.
Let $F$ be such a digraph and $F^*$ be a subdivision of $F$ in which an arc has been subdivided at least $3\cdot \mader_{\dic}(F)+1$ times. Since $\mader_{\dic}(F)\geq |V(F)| \geq 3$, the following result implies that the set of $(\mader_{\dic}(F))$-dicritical digraphs that do not contain any subdivision of $F^*$ is infinite.

\begin{restatable}{theorem}{propkdicritiquewithnolongdirectedpath}
    \label{thm:dicritical_nolong_directed_path}
    For every integer $k \geq 3$, there are infinitely many $k$-dicritical digraphs without any directed path on $3k+1$ vertices.
\end{restatable}

Theorem~\ref{thm:dicritical_nolong_directed_path} establishes a distinction between the directed and undirected cases. In the undirected case, for every fixed $k\geq 3$, there exists a non-decreasing function $f_k\colon \mathbb{N} \xrightarrow{} \mathbb{N}$ such that every $k$-critical graph on at least $f_k(\ell)$ vertices contains a path on $\ell$ vertices. This was first proved by Kelly and Kelly~\cite{kellyAJM76} in 1954, answering a question of Dirac. The bound on $f_k$ was then improved by Alon, Krivelevich and Seymour~\cite{alonJGT35} and finally settled by Shapira and Thomas~\cite{shapiraAM227}, who proved that the largest cycle in a $k$-critical graph on $n$ vertices has length at least $c_k \cdot \log(n)$, where $c_k$ is a constant depending only on $k$. This bound is best possible up to the multiplicative constant $c_k$, as shown by a construction of Gallai~\cite{gallai63a, gallai63b} (see~\cite{shapiraAM227}).

On the positive side, we adapt the proof of Alon et al.~\cite{alonJGT35} and show that, in $k$-dicritical digraphs, the length of the longest oriented cycle (\textit{i.e.} the longest cycle in the underlying graph) grows with the number of vertices, and so does the length of its longest oriented path (\textit{i.e.} the longest path in the underlying graph).

\begin{restatable}{theorem}{thmkdicriticallongestorientedpath}
    \label{thm:dicritical_longest_oriented_path}
    For every fixed integers $k\geq 2$ and $\ell\geq 3$, there are finitely many $k$-dicritical digraphs with no oriented cycle on at least $\ell$ vertices.
\end{restatable}

Since the digraphs constructed in Theorem~\ref{thm:dicritical_nolong_directed_path} contains many directed triangles, we propose to restrict ourselves to digraphs with large digirth.
The \textit{girth} of a graph $G$, denoted by $\girth(G)$, is the length of a smallest cycle in $G$, with the convention $\girth(G)=+\infty$ if $G$ is a forest.
The girth of a digraph $D$, denoted by $\girth(D)$, is the girth of its underlying graph.
The \textit{digirth} of a digraph is the length of its shortest directed cycle. By convention we have $\digirth(D)=+\infty$ if $D$ is acyclic.
A celebrated result of Erd\H{o}s (see~\cite{alon2016book}) states that there exist graphs of arbitrarily large chromatic number and arbitrarily large girth. It has been generalised by Bokal et al.~\cite{bokal2004circular} who showed the existence of digraphs of arbitrarily large digirth and dichromatic number.
For every integer $g$, we denote by $\mader_{\dic}^{(g)}(F)$ the least integer $k$ such that every digraph $D$ satisfying $\dic(D) \geq k$ and $\digirth(D) \geq g$ contains a subdivision of $F$. Note that $\mader_{\dic}(F) = \mader_{\dic}^{(2)}(F)$ and that $\mader_{\dic}^{(g)}(F)$ is non-increasing in $g$. 

\begin{conjecture}
\label{conj:mader_digirth}
For every digraph $F$ and every subdivision $F^*$ of $F$, there exists $g$ such that $\mader_{\dic}^{(g)}(F^*) \leq \mader_{\dic}(F)$.
\end{conjecture}

In order to provide some support to this conjecture, in Section~\ref{section:mader_bounded_by_F} we show that the value of $\mader_{\dic}^{(g)}(F^*)$ depends only on $F$ when $g$ is large enough. Our proof is strongly based on the key-lemma of~\cite{aboulkerEJC26}.

\begin{restatable}{theorem}{thmmadersubdivision}
    \label{thm:mader_subdivision}
    Let $k \geq 1$ be an integer. For every non-empty digraph $F$, if $F^*$ is obtained from $F$ by subdividing every arc at most $k-1$ times, then $\mader_{\dic}^{(k)}(F^*) \leq \frac{1}{3}\left(4^{m(F)+1}n(F)-1\right)$.
\end{restatable}

In Section~\ref{section:mader_bid_trees}, we prove that Conjecture~\ref{conj:mader_digirth} holds for every digraph $F$ whose underlying undirected graph $\UG(F)$ is a forest.

\begin{restatable}{theorem}{thmmaderbidtrees}
\label{thm:mader_digirth_bid_trees}
Let $k \geq 1$ be an integer and let $T$ be a bidirected tree. If $T^*$ is obtained from $T$ by subdividing every arc at most $k-1$ times, then $\mader_{\dic}^{(2k)}(T^*) \leq \mader_{\dic}(T) = n(T)$.
\end{restatable}

In the case of $T$ being an oriented tree, we improve Theorem~\ref{thm:mader_digirth_bid_trees} by proving $\mader_{\dic}^{(k)}(T) \leq \mader_{\dic}(T)$.
Observe that every digraph $D$ contains a subdigraph $H$ such that $\delta^+(H) \geq \dic(D) - 1$ (by taking $H$ a $\dic(D)$-dicritical subdigraph of $D$). 
Hence, for every integer $k$, if a digraph $F$ is such that every digraph $D$ with $\delta^+(D) \geq k-1$ contains a subdivision of $F$, then $\mader_{\dic}(F) \leq k$.
In Section~\ref{section:mader_delta}, we look for similar results using $\delta^+$ instead of $\dic$. 
Conjecture~\ref{conj:mader_digirth} for $F=C(1,2)$ appears to be a consequence of the following theorem (recall that $\mader_{\dic}(C(1,2)) = 3$).

\begin{restatable}{theorem}{thmspindleoutdegree}
\label{thm:mader_spindle_d+}
Let $k$ be an integer with $k \geq 2$. 
Every digraph $D$ with $\delta^+(D) \geq 2$ and $\digirth(D) \geq 8k-6$ contains a subdivision of $C(k,k)$.
\end{restatable}

When $k=2$, we improve Theorem~\ref{thm:mader_spindle_d+} by showing that every digraph $D$ with $\delta^+(D) \geq 2$ and $\digirth(D)\geq 3$ contains a subdivision of $C(2,2)$.

We finally consider out-stars. For integers $k,\ell$, let $S_k^{+(\ell)}$ be the digraph consisting of $k$ directed paths of length $\ell$ sharing their origin (and no other common vertices). The \textit{centre} of $S_k^{+(\ell)}$ is its unique source.

\begin{restatable}{theorem}{thmoutstar}
\label{thm:out_star}
Let $k$ and $\ell$ be two integers with $k\geq 2$ and $\ell\geq 1$. Every digraph $D$ with $\delta^+(D) \geq k$ and $\digirth(D) \geq \frac{k^\ell - 1}{k - 1}+1$ contains a copy of $S_k^{+(\ell)}$ with centre $u$ for every chosen vertex $u$.
\end{restatable}

When $k=2$, we show that Theorem~\ref{thm:out_star} can be improved by reducing the bound on $\digirth(D)$ down to $2\ell$.
We conclude in Section~\ref{sec:conclusion} by some open problems and further research directions.

\subsection*{Notations}

We refer the reader to~\cite{bang2009} for notation and terminology not explicitly defined in this paper.

Given a digraph $D$ (resp. a graph $G$), we denote by $V(D)$ (resp. $V(G)$) its set of vertices and $A(D)$ (resp. $E(G)$) its set of arcs (resp. set of edges). The order of $D$ is denoted by $n(D)=|V(D)|$ and its number of arcs is denoted by $m(D)=|A(D)|$.
For every vertex $v \in V(D)$, $N^+(v),N^-(v)$ denote respectively the out-neighbourhood and the in-neighbourhood of $v$ in $D$, $d^+(v)=|N^+(v)|, d^-(v)=|N^-(v)|$, and $d(v)=d^+(v)+d^-(v)$.
Moreover, for every set $A$ of vertices, we write $N^+(A) = \bigcup_{a \in A} N^+(a)$, $N^-(A) = \bigcup_{a \in A} N^-(a)$ and $N(A) = N^+(A) \cup N^-(A)$.
We also define $\delta^+(D) = \min_{v \in V(D)} d^+(v)$, $\delta^-(D) = \min_{v \in V(D)} d^-(v)$, and $\delta(D) = \min_{v \in V(D)} d(v)$. We similarly define $\Delta(D) = \max_{v \in V(D)} d(v)$.

A {\it digon} in $D$ is a pair of opposite arcs between two vertices. Such a pair of arcs $\{uv,vu\}$ is denoted by $[u,v]$.
The {\it underlying graph} of a digraph $D$, denoted by $\UG(D)$, is the undirected graph on the same vertex set which contains an edge linking two vertices if $D$ contains at least one arc linking these two vertices.
An {\it oriented graph} is a digraph with digirth at least $3$.
A {\it bidirected digraph} is any digraph $D$ such that for every $x,y \in V(D)$, $xy \in A(D) \Leftrightarrow yx \in A(D)$,
and we write $D = \bid{G}$ where $G = \UG(D)$.
A \textit{cycle} is an undirected connected graph in which every vertex has degree $2$.
For such a cycle $C$ with vertices $u_1, \dots, u_\ell$, we write $C = (u_1, \dots, u_\ell,u_1)$ if $N(u_i) = \{u_{(i-1) \mod \ell},u_{(i+1) \mod \ell}\}$ for every $i \in [\ell]$.
An \textit{oriented cycle} (resp. \textit{oriented path}) is an oriented graph whose underlying graph is a cycle (resp. path).
A \textit{directed cycle} is an oriented cycle in which every vertex has in- and out-degree $1$.
For such a directed cycle $C$ with vertices $u_1, \dots, u_\ell$, we write $C = (u_1, \dots, u_\ell,u_1)$ if $N^+(u_i) = \{u_{(i+1) \mod \ell}\}$ for every $i \in [\ell]$.
The directed cycle of length $3$ is also called the {\it directed triangle}.
A \textit{directed path} is obtained from a directed cycle by removing a vertex. The directed cycle on $n$ vertices and the directed path on $n$ vertices are respectively denoted by $\ori{C_n}$ and $\ori{P_n}$. The {\it lengths} of $\ori{P_n}$ and $\ori{C_n}$ are respectively $n-1$ and $n$.
For such a directed cycle $P$ with vertices $u_1, \dots, u_\ell$, we write $P = (u_1, \dots, u_\ell)$ if $N^+(u_i) = \{u_{i+1}\}$ for every $i \in [\ell-1]$.
An \textit{antidirected path} is an orientation of a path where every vertex $v$ satisfies $\min(d^+(v),d^-(v)) = 0$.

A {\it subdivision} of a digraph $F$ is any digraph obtained from $F$ by replacing every arc $uv$ by a directed path from $u$ to $v$.
If a digraph $D$ contains a subdivision of $F$, we say that {\it $D$ contains $F$ as a subdivision}.

If $u$ is a vertex of a digraph $D$ with $d^-(u)=1$, we define the {\it predecessor} of $u$ in $D$, denoted by $\pred_D(u)$, as the unique in-neighbour of $u$ in $D$.
Similarly, if $d^+(u)=1$ , we define the {\it successor} of $u$ in $D$, denoted by $\succ_D(u)$, as the unique out-neighbour of $u$ in $D$.

Given two vertices $a,b$ in a directed cycle $C$ (with possibly $a=b$), we denote by $C[a,b]$ the directed path from $a$ to $b$ along $C$ (which is the single-vertex path when $a=b$).
Moreover, we define $C[a,b[ = C[a,b]-\{b\}$, $C]a,b] = C[a,b]-\{a\}$, and $C]a,b[ = C[a,b]-\{a,b\}$. 
Note that these subpaths may be empty.
Given a directed path $P$ and two vertices $a,b$ in $V(P)$, we use similar notations $P[a,b],P[a,b[,P]a,b]$ and $P]a,b[$.
We denote by $\init(P)$ the first vertex of $P$ (\textit{i.e.} the unique vertex with in-degree $0$) and $\term(P)$ its last one.
Given two directed paths $P,Q$ such that $V(P) \cap V(Q) = \{x\}$ where $x = \term(P) = \init(Q)$, the {\it concatenation} of $P$ and $Q$, denoted by $P\cdot Q$, is the digraph $(V(P) \cup V(Q), A(P) \cup A(Q))$.
The vertices in $V(P) \setminus \{\init(P),\term(P)\}$ are called the {\it internal vertices} of $P$.
If $U$ and $V$ are two  sets of vertices in $D$, then a {\it $(U,V)$-path} in $D$ is a directed path $P$ in $D$ with $\init(P)\in U$ and $\term(P) \in V$, and we also say that $P$ is a directed path from $U$ to $V$. If $U=\{u\}$ (resp. $V=\{v\}$) then we simply write $u$ for $U$ (resp. $v$ for $V$) in these notations.
The {\it distance} from $u$ to $v$, denoted by $\dist(u,v)$, is the length of a shortest $(u,v)$-path, with the convention $\dist(u,v) = +\infty$ if no such path exists.

A digraph is {\it connected} if its underlying graph is connected. It is {\it strongly connected} if for every ordered pair $(u,v)$ of its vertices, there exists a directed path from $u$ to $v$.

\section{An improved bound on \texorpdfstring{$\mader_{\dic}(\bid{K_n})$}{mader(Kn)}}
\label{section:mader_Kn}

This section is devoted to the proof of Proposition~\ref{thm:mader_Kn}. We need the two following lemmas.

\begin{lemma}[Aboulker et al.~{\cite[Lemma~31]{aboulkerEJC26}}]
\label{lemma:aboulker_et_al_add_arc}
For every integer $k$ and every digraph $D$ with $\dic(D) \geq 4k-3$, there is a subdigraph $H$ of $D$ with $\dic(H) \geq k$ such that
for every pair $u,v$ of distinct vertices in $H$, there is a directed path from $u$ to $v$ in $D$ whose internal vertices are in $V(D) \setminus V(H)$.

In particular, for every digraph $F$ and every arc $e$ in $F$,
\[\mader_{\dic}(F) \leq 4\cdot \mader_{\dic}(F\setminus e)-3.\]
\end{lemma}
We skip the proof of the following easy lemma.

\begin{lemma}\label{lemma:disjoint_union}
If $F_1 + F_2$ denotes the disjoint union of two digraphs $F_1$ and $F_2$, then
$\mader_{\dic}(F_1+F_2) \leq \mader_{\dic}(F_1) + \mader_{\dic}(F_2)$.
\end{lemma}

We are now ready to prove Proposition~\ref{thm:mader_Kn}, let us first restate it.

\begin{proposition}\label{thm:mader_Kn}
$\mader_{\dic}(\bid{K_n}) \leq 4^{\frac{2}{3}n^2 + 2n - \frac{8}{3}}$.
\end{proposition}

\begin{proof}
Let $f(n) = \mader_{\dic}(\bid{K_n})$ for every $n \geq 1$.
Clearly $f(1)=1$.
Let $g(x) = 4^{\frac{2}{3}x^2+2x - \frac{8}{3}}$ for every positive real $x$.
Observe that $g$ is non-decreasing.
We will show by induction on $n$ that $f(n) \leq g(n)$ for every positive integer $n$.
For $n=1$, $f(1) = 1 = g(1)$.
Now suppose $n \geq 2$.

If $n$ is even, by Lemmas~\ref{lemma:aboulker_et_al_add_arc} and~\ref{lemma:disjoint_union} we deduce the following inequalities.
\begin{align*}
f(n) &\leq 4^{\frac{n^2}{2}} \cdot \mader_{\dic}(\bid{K_{\frac{n}{2}}} + \bid{K_{\frac{n}{2}}}) \\
&\leq 4^{\frac{n^2}{2}} \cdot 2 \cdot f\left(\frac{n}{2}\right) \\
&\leq 4^{\frac{n^2}{2}} \cdot 2 \cdot g\left(\frac{n}{2}\right) \\
&\leq 4^{\frac{n^2}{2} + n} \cdot g\left(\frac{n}{2}\right). \\
\end{align*}
If $n$ is odd, then $f(n) \leq 4^{2(n-1)}\mader_{\dic}(\bid{K_1}+\bid{K_{n-1}}) \leq 4^{2(n-1)}(1+f(n-1)) \leq 4^{2n-1}f(n-1)$.
Hence:
\begin{align*}
f(n) &\leq 4^{2n-1}f(n-1) \\
&\leq 4^{2n-1} \cdot 4^{\frac{(n-1)^2}{2}} \cdot \mader_{\dic}(\bid{K_{\frac{n-1}{2}}} + \bid{K_{\frac{n-1}{2}}})  \\
&\leq 4^{2n-1} \cdot 4^{\frac{(n-1)^2}{2}} \cdot 2 \cdot f\left(\frac{n-1}{2}\right) \\
&\leq 4^{2n-1} \cdot 4^{\frac{(n-1)^2}{2}} \cdot 2 \cdot g\left(\frac{n-1}{2}\right) \\
&\leq 4^{\frac{n^2}{2} + n - \frac{1}{2}} \cdot 2 \cdot g\left(\frac{n}{2}\right) \\
&= 4^{\frac{n^2}{2} + n} \cdot g\left(\frac{n}{2}\right). \\
\end{align*}
In both cases we have:
\[
f(n) \leq 4^{\frac{n^2}{2} + n} \cdot g\left(\frac{n}{2}\right) = 4^{n^2(\frac{1}{2}+\frac{1}{6}) + n(1 + 1)-\frac{8}{3}} = g(n).
\]
\end{proof}

\section{Paths and cycles in large dicritical digraphs}
\label{section:large_dicritical}

This section is devoted to the proofs of Theorems~\ref{thm:dicritical_nolong_directed_path} and~\ref{thm:dicritical_longest_oriented_path}. We first prove the following useful observation.

\begin{lemma}\label{lemma:universalvertex}
    For every integer $k$, if $D$ is a $k$-dicritical digraph and if $D'$ is obtained from $D$ by adding a vertex $u$ with $N^+(u)=N^-(u)=V(D)$,
    then $D'$ is $(k+1)$-dicritical.
\end{lemma}

\begin{proof}
First we show that $\chi(D') \geq k+1$.
Indeed, if $\phi\colon V(D') \to [k]$ is a $k$-dicolouring of $D'$, then $\phi(v) \neq \phi(u)$ for every $v \in V(D)$, and so $\phi$ induces
a $(k-1)$-dicolouring of $D$, a contradiction.

It remains to show that for every arc $vw$ in $D'$, $D'\setminus vw$ is $k$-dicolourable.
If $vw \in A(D)$, then since $D$ is $k$-dicritical, $D \setminus vw$ admits a $(k-1)$-dicolouring $\phi\colon V(D) \to [k-1]$.
Then extending $\phi$ to $V(D')$ by $\phi(u)=k$ yields a $k$-dicolouring of $D'$.
If $u\in \{v,w\}$, then consider a $k$-dicolouring $\phi\colon V(D)\setminus \{v,w\} \to [k-1]$ of $D' - \{v,w\}$, which exists since $D' - \{v,w\}$ is a proper subdigraph of $D$.
Now set $\phi(v)=\phi(w)=k$. Colour $k$ induces an acyclic digraph of $D'\setminus vw$, and this yields a $k$-dicolouring of $D'\setminus vw$.
\end{proof}

\propkdicritiquewithnolongdirectedpath*

\begin{proof}
    Let $n$ be an odd integer.
    Let $D_{k,n}$ be the digraph constructed as follows.
    Start with the antidirected path $P=(p_1, \dots, p_n)$ on $n$ vertices in which $d^+(p_1) = 1$.
    Add the digon $[p_1,p_n]$, and two vertices $x_1, x_2$ with a digon $[x_1,x_2]$.
    For every arc $uv$ of $P$, add the arcs $v x_i,x_i u$ for every $i \in [2]$.
    Finally, add $k-3$ vertices $x_3, \dots, x_{k-1}$ inducing a copy of $\bid{K_{k-3}}$ and add the digon $[x_{i+2},u]$ for every $i \in [k-3]$ and every $u \in V(P) \cup \{x_1,x_2\}$.
    See Figure~\ref{fig:D_n_bis} for an illustration.

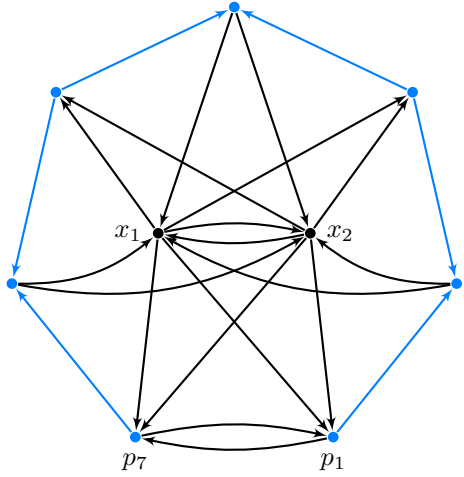
\begin{figure}[hbtp]
    \begin{center}	
          \begin{tikzpicture}[thick,scale=1, every node/.style={transform shape}]
            \tikzset{vertex/.style = {circle,fill=black,minimum size=4pt, inner sep=0pt}}
            \tikzset{edge/.style = {->,> = latex'}}

            \node[vertex, g-blue, label=below:$p_1$] (u1) at (90-1080/7:3) {};
            \node[vertex, g-blue, label=below:$p_7$] (u7) at (2160/7 -1080/7 + 90:3) {};
            \foreach \i in {2,...,6}{
                \pgfmathtruncatemacro{\j}{\i-1}
                \node[vertex, g-blue] (u\i) at (\j*360/7 - 1080/7 + 90:3) {};
            }
            \node[vertex, label=left:$x_1$] (x) at (-1,0) {};
            \node[vertex, label=right:$x_2$] (y) at (1,0) {};
            \draw[edge, g-blue] (u1) to (u2);
            \draw[edge, g-blue] (u3) to (u2);
            \draw[edge, g-blue] (u3) to (u4);
            \draw[edge, g-blue] (u5) to (u4);
            \draw[edge, g-blue] (u5) to (u6);
            \draw[edge, g-blue] (u7) to (u6);
            \draw[edge] (x) to (u1); \draw[edge] (y) to (u1);
            \draw[edge] (x) to (u3); \draw[edge] (y) to (u3);
            \draw[edge] (x) to (u5); \draw[edge] (y) to (u5);
            \draw[edge] (x) to (u7); \draw[edge] (y) to (u7);
            \draw[edge, bend left=20] (u2) to (x); \draw[edge, bend left=20] (u2) to (y);
            \draw[edge] (u4) to (x); \draw[edge] (u4) to (y);
            \draw[edge, bend right=20] (u6) to (x); \draw[edge, bend right=20] (u6) to (y);
            \draw[edge, bend left=12] (u1) to (u7);
            \draw[edge, bend left=12] (u7) to (u1);
            \draw[edge, bend left=12] (x) to (y);
            \draw[edge, bend left=12] (y) to (x);
          \end{tikzpicture}
      \caption{The digraph $D_{3,7}$. The antidirected path $P$ is in blue.}
      \label{fig:D_n_bis}
    \end{center}
\end{figure}

    Let us show that $D_{k,n}$ is $k$-dicritical.
    Since $D_{k,n}$ is obtained from $D_{3,n}$ by adding $k-3$ vertices linked to all other vertices by a digon, by Lemma~\ref{lemma:universalvertex}
    it is enough to show that $D_{3,n}$ is $3$-dicritical.
    
    First we show that $\dic(D_{3,n}) > 2$.
    Suppose for contradiction that there is a $2$-dicolouring $\phi\colon V(D_{3,n}) \to [2]$ of $D_{3,n}$.
    Without loss of generality, $\phi(x_1)=1$ and $\phi(x_2)=2$.
    For every arc $uv$ of $P$, for every $i \in [2]$, $uvx_iu$ is a directed triangle. This implies that $\phi(u) \neq \phi(v)$.
    Since $P$ has an odd number of vertices, $\phi(p_1)=\phi(p_n)$, which is a contradiction as $[p_1,p_n]$ is a digon in $D_{3,n}$.

    Let $uv$ be an arc in $D_{3,n}$. We show that $\dic(D_{3,n} \setminus uv) \leq 2$.
    If $\{u,v\} = \{x_1,x_2\}$, set $\phi(x_1)=\phi(x_2)=\phi(p_1)=1$, $\phi(p_2)= \dots = \phi(p_{n})=2$.
    If $\{u,v\} = \{p_1,p_n\}$, set $\phi(x_1)=1$, $\phi(x_2)=2$, $\phi(p_i) = 1$ if $i$ is even, $\phi(p_i)=2$ if $i$ is odd.
    If $uv \in A(P)$, set $\phi(x_1)=1$, $\phi(x_2)=2$, and colour $V(P)$ such that the only monochromatic pair of adjacent vertices in $P$ is $\{u,v\}$.
    If $u \in V(P)$ and $v = x_i$, set $\phi(x_1)=1$, $\phi(x_2)=2$, $\phi(u)=i$, and colour $V(P-u)$ such that two adjacent vertices in $V(P-u)$ receive distinct colours.
    The other cases are symmetric.
    In each case, one can check that this gives a proper $2$-dicolouring of $D_{3,n}\setminus  uv$.

    It remains to prove that $D_{k,n}$ does not contain a directed path on $3k+1$ vertices.
    Let $Q$ be a directed path in $D_{k,n}$.
     Let $y_1, \dots, y_\ell$ be the vertices of $V(Q) \cap \{x_1, \dots, x_{k-1}\}$ in order of appearance along $Q$. 
    Let $Q_i$ be the subpath $Q]y_i,y_{i+1}[$ of $Q$ for every $i \in[\ell-1]$ and let $Q_0 = Q[\init(Q),y_1[$ and $Q_\ell = Q]y_\ell,\term(Q)]$). 
    Note that some $Q_i$s may be empty.
    
    Then observe that at most one of the $Q_i$s intersects both $p_1$ and $p_n$.
    Except this one, which has at most three vertices, all the $Q_j$s have at most two vertices (because $P$ is anti-directed).
    We conclude that the number of vertices in $P$ is at most $3\ell+3 \leq 3k$.
\end{proof}

\thmkdicriticallongestorientedpath*
\begin{proof}
    Let $k\geq 2$ be a fixed integer. We will show the existence of a function $f_k\colon \mathbb{N} \xrightarrow[]{} \mathbb{N}$ such that every $k$-dicritical digraph 
 on at least $f_k(\ell)$ vertices contains an oriented path on $\ell$ vertices. We will then use a result of Dirac to show that every $k$-dicritical digraph on at least $f_k(\frac{1}{4}\ell^2)$ vertices contains an oriented cycle on $\ell$ vertices, implying the result.

    Given a digraph $H$, $\cc(H)$ is the number of connected components of $H$ (\textit{i.e.} the number of connected components of $\UG(H)$). Our proof is strongly based on the following claim. 
    \addtocounter{theorem}{-9}
    \begin{claim}
        \label{claim:number_connected_component_dicritical}
        Let $D=(V,A)$ be a $k$-dicritical digraph and $S\subseteq V$, then $\cc(D-S) \leq (k-1)^{|S|} \cdot 3^{\binom{|S|}{2}}$.
    \end{claim}
    \addtocounter{theorem}{9}
    \begin{proofclaim}
        Assume this is not the case, \textit{i.e.} there exists a $k$-dicritical digraph $D$ and a subset of its vertices $S$ such that $\cc(D-S) >(k-1)^{|S|} \cdot 3^{\binom{|S|}{2}}$. We denote by $H_1,\dots,H_r$ the connected components of $D-S$. 
        
        For every $i\in [r]$, let $\alpha_i$ be a $(k-1)$-dicolouring of $D - V(H_i)$, the existence of which is guaranteed by the dicriticality of $D$. 
        Let $s=|S|$ and $v_1,\dots,v_s$ be any fixed ordering of $S$. For every $i\in[r]$, we let $\sigma_i^1$ be the ordered set $(\alpha_i(v_1),\dots,\alpha_i(v_s))$.
        We also define $\sigma_i^2$ as the set of all ordered pairs $(u,v)\in S^2$ such that $D-V(H_i)$, coloured with $\alpha_i$, contains a monochromatic directed path from $u$ to $v$. We finally define the $i^{\text{th}}$ configuration $\sigma_i$ as the ordered pair $(\sigma_i^1, \sigma_i^2)$.
        
        For every pair of vertices $u,v$ in $S$ and every $i\in[r]$, note that at most one of the ordered pairs $(u,v),(v,u)$ actually belongs to $\sigma_i^2$, for otherwise $D- V(H_i)$, coloured with $\alpha_i$, contains a monochromatic directed cycle. Hence, the number of distinct configurations is at most $(k-1)^{s} \cdot 3^{\binom{s}{2}}$.
        By the pigeonhole principle, since $r > (k-1)^{s} \cdot 3^{\binom{s}{2}}$, there exist two distinct integers $i,j\in [r]$ such that $\sigma_i = \sigma_j$. Let $\alpha$ be the colouring of $D$ defined as follows:
        \[
        \alpha(v) = \left\{
            \begin{array}{ll}
                \alpha_i(v) & \mbox{if } v \in (V(D) \setminus V(H_i)) \\
                \alpha_j(v) & \mbox{otherwise.}
            \end{array}
        \right.
        \]

        We claim that $\alpha$ is a $(k-1)$-dicolouring of $D$. Assume for a contradiction that it is not, so $D$, coloured with $\alpha$, contains a monochromatic directed cycle. Among all such cycles $C$, we choose one for which the size of $V(C) \cap V(H_i)$ is minimised. If $V(C) \cap V(H_i) = \emptyset$, then $C$ is a monochromatic directed cycle of $D - H_i$, a contradiction to the choice of $\alpha_i$. Analogously, we have $V(C) \setminus V(H_i) \neq \emptyset$ by choice of $\alpha_j$. 
        
        Assume first that $|V(C) \setminus V(H_i)| = 1$, implying that $C$ contains exactly one vertex $s$ in $S$ and $(V(C) \setminus \{s\}) \subseteq V(H_i)$. Since $\sigma_i^1 = \sigma_j^1$, we have $\alpha_i(s) = \alpha_j(s)$, which implies that $C$ is a monochromatic directed cycle of $D - H_j$ coloured with $\alpha_j$, a contradiction to the choice of $\alpha_j$.

        Henceforth we can assume that $C$ contains a directed path $P$ on at least three vertices, with initial vertex $u$ and terminal vertex $v$, such that $V(P) \cap S =\{u,v\}$ and $V(P) \subseteq (V(H_i) \cup \{u,v\})$. The existence of $P$ ensures that $(u,v)$ belongs to $\sigma_j^2$. Hence, since $\sigma_i = \sigma_j$, there exists a monochromatic directed path $P'$ in $D-V(H_i)$ coloured with $\alpha_i$, from $u$ to $v$, and with the same color as $P$. Hence, replacing $P$ by $P'$ in $C$, we obtain a closed walk which, coloured with $\alpha$, contains a monochromatic directed cycle $C'$ such that $|V(C') \cap V(H_i)| < |V(C) \cap V(H_i)|$, a contradiction to the choice of $C$.
    \end{proofclaim}

    We are now ready to prove the existence of $f_k$. Let $D$ be a $k$-dicritical digraph whose underlying graph $G$ does not contain any path on $\ell$ vertices. Let $v$ be any vertex of $D$ and $T$ be a spanning DFS-tree of $G$ rooted in $v$ (recall that $D$ is connected since it is dicritical). Let $h$ be the depth of $T$ (\textit{i.e.} the maximum number of vertices in a branch of $T$), then $h$ is at most $\ell$ since $G$ does not contain any path of length $\ell$. 
    
    For every vertex $x$, let $S_x$ be the ancestors of $x$ (including $x$ itself) in $T$ and $d_T(x)$ be the number of children of $x$ in $T$. Since $T$ is a DFS-tree, note that for every neighbour $y$ of $x$, $x$ and $y$ must belong to the same branch. Hence, $d_T(x) \leq \cc(D-S_x)$. Since $|S_x| \leq h \leq \ell$, we deduce from Claim~\ref{claim:number_connected_component_dicritical} that $d_T(x) \leq (k-1)^{\ell} \cdot 3^{\binom{\ell}{2}}$.
    Since $T$ is spanning, we obtain that $|V(D)| \leq \left( (k-1)^{\ell} \cdot 3^{\binom{\ell}{2}} \right)^{\ell-1} = f_k(\ell) - 1$.

    \medskip

    Dirac proved that every 2-connected graph that contains a path of length $t$ actually contains a cycle of length at least $2\sqrt{t}$ (see~\cite[Problem~10.29]{lovaszbook2007}). It is straightforward to show that every $k$-dicritical digraph is $2$-connected. Hence, if $D$ is a $k$-dicritical digraph on at least $f_k(\frac{1}{4}\ell^2)$ vertices, then $D$ contains an oriented cycle of length at least $\ell$, implying the result.
\end{proof}

\section{Subdivisions in digraphs with large digirth}
\label{section:mader_bounded_by_F}

This section is devoted to the proof of Theorem~\ref{thm:mader_subdivision}.

\thmmadersubdivision*

\begin{proof}
    We proceed by induction on $m(F)$, the result being trivial when $m(F)=0$. Let $F$ be any digraph with $m>0$ arcs, and let $F^*$ be a digraph obtained from $F$ by subdividing every arc at most $k-1$ times.
    
    Let $uv\in A(F)$ be any arc, and $P=x_1,\dots,x_r$ its corresponding directed path in $F^*$ (where $u=x_1$ and $v=x_r$). Then we only have to prove that $\mader_{\dic}^{(k)}( F^* \setminus x_1x_2 ) \leq \frac{4^{m(F)}n(F)-1}{3} + 1$. If this is true, then by Lemma~\ref{lemma:aboulker_et_al_add_arc}, we get that  $\mader_{\dic}^{(k)}( F^* ) \leq 4\left(\frac{4^{m(F)}n(F)-1}{3} + 1\right) - 3$ which shows the result.
    
    Let $D$ be any digraph with dichromatic number at least $\frac{4^{m(F)}n(F)-1}{3} + 1$ and digirth at least $k$ and let $B\subseteq V(D)$ be a maximal acyclic set in $D$. Then $\dic(D - B) \geq \frac{4^{m(F)}n(F)-1}{3}$, so by induction $D-B$ must contain a subdivision of $F\setminus uv$ where each arc has been subdivided at least $k-1$ times. This is also a subdivision of $F^* - \{ x_2,\dots,x_{r-1} \}$. Let $y$ be the vertex in $D-B$ corresponding to $x_r$. By maximality of $B$, there must be a directed cycle $C$ in $D$ such that $V(C) \cap V(D-B) = \{y\}$. Note that $C$ has length at least $k$. Thus, ignoring the leaving arc of $y$ in $C$, we have found a subdivision of $F^*\setminus x_1x_2$ in $D$, showing the result. 
\end{proof}

\section{Subdivisions of trees in digraphs with large digirth}
\label{section:mader_bid_trees}

This section is devoted to the proofs of Theorems~\ref{thm:mader_digirth_bid_trees}.

\thmmaderbidtrees*

\begin{proof}
    We proceed by induction on $n(T)$.
    Suppose $n(T) \geq 2$, the result being trivial when $n(T) = 1$.
    Let $f$ be a leaf of $T$ with neighbour $p$, and we denote by $(T-f)^*$
    the bidirected tree $T-f$ with every arc subdivided exactly $k-1$ times.
    By induction hypothesis $\mader_{\dic}^{(2k)}((T-f)^*) \leq \mader_{\dic}(T-f) = n(T)-1$.
    Let $D$ be a digraph with $\digirth(D) \geq 2k$ and $\dic(D) \geq n(T)$, and consider a maximal acyclic set $A$ in $D$.
    Then $\dic(D-A) \geq n(T)-1$ and so by induction hypothesis, $D-A$ contains a subdivision of $(T-f)^*$.
    Let $y \in V(D)\setminus A$ be the vertex corresponding to $p \in V(T)$ in the subdivision of $(T-f)^*$ contained in $D-A$.
    By maximality of $A$, $A+y$ contains a directed cycle $C$ with $V(C) \setminus A = \{y\}$.
    As $\digirth(D) \geq 2k$, $C$ has length at least $2k$.
    Then the subdivision of $(T-f)^*$ in $D-A$ together with $C$ gives the desired subdivision
    of $T^*$.
\end{proof}

As mentioned in the introduction, if $T$ is an oriented tree, then Theorem~\ref{thm:mader_digirth_bid_trees} can be strengthened as follows.

\begin{theorem}\label{thm:univ_digirth_oriented_trees}
    Let $k \geq 1$ be an integer and let $T$ be an oriented tree. If $T^*$ is obtained from $T$ by subdividing every arc at most $k-1$ times, then
\end{theorem}

\begin{proof}
    We proceed by induction on $n(T)$. 
    Suppose $n(T) \geq 2$, the result being trivial when $n(T) = 1$.
    
    For every arc $e$ of $T$, we denote by $s(e)$ the number of subdivisions of $e$ in $T^*$. Let $f$ be a leaf of $T$ with neighbour $p$, and we denote by $(T-f)^*$
    the oriented tree $T-f$ where every arc $e$ is subdivided $s(e)\leq k-1$ times.
    
    Let $D$ be a digraph with $\digirth(D) \geq k$ and $\dic(D) \geq n(T)$, and consider a maximal acyclic set $A$ in $D$. We have $\dic(D-A) \geq n(T)-1$ and, by the induction hypothesis, $D-A$ contains a copy  of $(T-f)^*$.
    Let $y \in V(D)\setminus A$ be the vertex corresponding to $p \in V(T)$ in the copy of $(T-f)^*$ contained in $D-A$. By maximality of $A$, $A+y$ contains a directed cycle $C$ with $V(C) \setminus A = \{y\}$. As $\digirth(D) \geq k$, $C$ has length at least $k$. 
    
    If the arc between $p$ and $f$ goes from $p$ to $f$ then we define $P$ as the directed path on $s(pf)$ vertices starting from $y$ along $C$.
    Otherwise, it goes from $f$ to $p$ and then we define $P$ as the directed path on $s(fp)$ vertices, along $C$, ending on $y$.
    In both cases, the copy of $(T-f)^*$ in $D-A$ together with $P$ gives the desired copy of $T^*$.
\end{proof}

\section{Subdivisions in digraphs of large out-degree and large digirth}
\label{section:mader_delta}

\subsection{Subdivisions of \texorpdfstring{$C(k,k)$}{C(k,k)}}

This section is devoted to the proof of Theorem~\ref{thm:mader_spindle_d+}.
\thmspindleoutdegree*

\addtocounter{theorem}{-6}

\begin{proof}
We will prove the following stronger statement: for every digraph $D$ with $\digirth(D) \geq 8k-6$ and $v_0 \in V(D)$,
if $d^+(v_0) \geq 1$ and $d^+(v) \geq 2$ for every $v \in V(D) \setminus\{v_0\}$, then $D$
contains a subdivision of $C(k,k)$.
We now consider a counterexample to this statement with minimum number of vertices, and minimum number of arcs if equality holds.

\begin{claim}
    $D$ is strongly connected.
\end{claim}

\begin{proofclaim}
    Let $C$ be a terminal strongly connected component of $D$, that is a strongly connected component such that there is no arc going out of $C$. Then $C$ is also a counterexample, so by minimality of $D$ we have $D=C$, and $D$ is strongly connected.
\end{proofclaim}

\begin{claim}
    $d^+(v) = 2$ for every vertex $v \neq v_0$ and $d^+(v_0)=1$.
\end{claim}

\begin{proofclaim}
    If $v \neq v_0$ is a vertex with at least $3$ out-neighbours $w_1,w_2,w_3$, then
    $D\setminus v w_3$ is a smaller counterexample.
    Similarly, if $d^+(v_0) > 1$, then $v_0$ has at least two distinct out-neighbours $w_1,w_2$,
    and $D \setminus vw_2$ is a smaller counterexample.
\end{proofclaim}

Given two vertices $u,v$ of $D$, a {\it $(u,v)$-vertex-cut} is a vertex $x\in V(D) \setminus \{u,v\}$ which intersects every $(u,v)$-path of $D$.

\begin{claim}\label{claim:dist_uv}
    Let $u,v$ be two vertices in $D$.
    If $\dist(u,v) \leq 7k-6$, then there exists a $(v,u)$-vertex-cut.
\end{claim}

\begin{proofclaim}
    Suppose the contrary for contradiction.
    Then by Menger's theorem, there exist two internally vertex-disjoint $(v,u)$-paths
    $P_1$ and $P_2$.
    As $\digirth(D) \geq 8k-6$, both $P_1$ and $P_2$ have length at least $k$, and so
    $P_1 \cup P_2$ is a subdivision of $C(k,k)$ with source $v$ and sink $u$.
\end{proofclaim}

For every directed cycle $C$ in $D$, let $\rho(C)$ be the number of vertices in the
largest connected component of $D-V(C)$.
We say that $C$ is isometric if for every $u,v \in V(C)$, $C$ contains
a shortest $(u,v)$-path in $D$.
Clearly $D$ contains an isometric cycle (it is enough to take a minimum directed cycle), and we consider among them an isometric cycle $C$ which maximises $\rho(C)$.

Let $ab$ be an arc along $C$.
Let $c_1, \dots, c_\ell$ be the $(b,a)$-vertex-cuts in $D$. Observe that at least one such vertex-cut exists by Claim~\ref{claim:dist_uv}.
As $C$ contains a $(b,a)$-path, all these vertices belong to $V(C)$,
and we suppose that they appear in this order $c_1, \dots ,c_\ell$ along $C$ starting at $b$.
By convention we also define $c_0=b$ and $c_{\ell+1}=a$.

\begin{claim}\label{claim:dist_ci_cj}
    $\dist(c_i,c_{i+1}) \leq k-1$ for every $i=0, \dots, \ell$.
\end{claim}

\begin{proofclaim}
    Suppose for contradiction that $\dist(c_i,c_{i+1}) \geq k$.
    Assume first that there exists a $(c_i,c_{i+1})$-vertex-cut $x$.
    We claim that $x$  is also a $(b,a)$-vertex-cut.
    Consider a $(b,a)$-path $Q$. Then $Q$ passes through $c_i$ and $c_{i+1}$ in this order, otherwise the concatenation of $Q[b,c_{i+1}]$ and $C[c_{i+1},a]$ is a $(b,a)$-path avoiding $c_i$, a contradiction. By definition of $x$, it belongs to $Q[c_i,c_{i+1}]$. Hence $x$ intersects every $(b,a)$-path, and so $x \in \{c_1, \dots, c_\ell\}$, a contradiction since $c_1, \dots, c_\ell$ are in this order along $C$.

    This shows, by Menger's theorem, that there are two internally vertex-disjoint $(c_i,c_{i+1})$-paths $P_1,P_2$.
    Then $P_1$ and $P_2$ have length at least $k$, and
    so $P_1 \cup P_2$ is a subdivision of $C(k,k)$ with source $c_i$ and sink $c_{i+1}$, a contradiction.
\end{proofclaim}

Let $i_0$ be the least index $i$ such that $\dist(b,c_i) \geq k$, and let $i_1$ be the largest index $i$ such that $\dist(c_i,a) \geq k$.
By choice of $i_0$, we have $\dist(b,c_{i_0-1}) \leq k-1$. 
By Claim~\ref{claim:dist_ci_cj}, we have $\dist(c_{i_0-1},c_{i_0}) \leq k-1$, which implies $\dist(b,c_{i_0})\leq \dist(b,c_{i_0-1}) + \dist(c_{i_0-1},c_{i_0}) \leq 2k-2$. Similarly we have $\dist(c_{i_1},a) \leq 2k-2$.
Therefore, we have 
\[
\dist(c_{i_0},c_{i_1}) = |V(C)| - \dist(c_{i_1},a) - 1 - \dist(b,c_{i_0}) \geq 4k-3,
\]
which implies $i_1 - i_0 \geq 5$ by Claim~\ref{claim:dist_ci_cj}. See Figure~\ref{fig:digirth_large_C_k_k} for an illustration. 

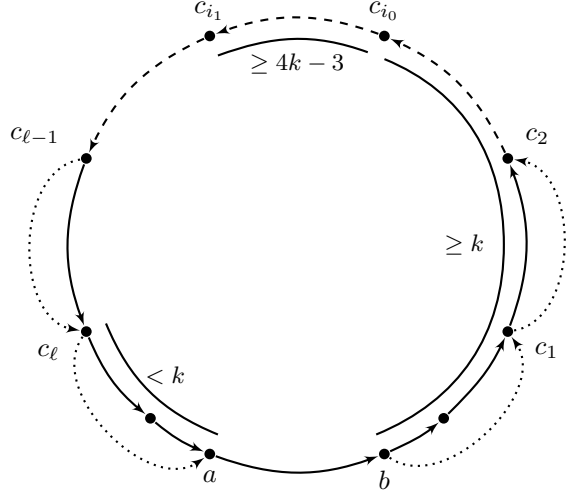
\begin{figure}[hbtp]
        \begin{center}	
              \begin{tikzpicture}[thick,scale=1, every node/.style={transform shape}]
                \tikzset{vertex/.style = {circle,fill=black,minimum size=4pt, inner sep=0pt}}
                \tikzset{edge/.style = {->,> = latex'}}
                
                \node[vertex, label=below:$a$] (v1) at (-112.5:3) {};
                \node[vertex, label=below:$b$] (v2) at (-67.5:3) {};
                \node[vertex] (v3) at (-50:3) {};
                \node[vertex] (vl) at (-130:3) {};
                \node[vertex, label={[label distance=0.4em]355:$c_1$}] (s1) at (-22.5:3) {};
                \node[vertex, label={[label distance=0.1em]30:$c_2$}] (s2) at (22.5:3) {};
                \node[vertex, label=above:$c_{i_0}$] (s3) at (67.5:3) {};
                \node[vertex, label=above:$c_{i_1}$] (s4) at (112.5:3) {};
                \node[vertex, label={[label distance=0.4em]170:$c_{\ell-1}$}] (srm1) at (157.5:3) {};
                \node[vertex, label={[label distance=0.4em]190:$c_\ell$}] (sr) at (-157.5:3) {};

                \draw[edge] (v1) to[out=-20, in=-160] (v2);
                \draw[edge] (v2) to[out=25, in=-140] (v3);
                \draw[edge] (v3) to[out=45, in=-115] (s1);
                \draw[edge] (s1) to[out=70, in=-70] (s2);
                \draw[edge,dashed] (s2) to[out=115, in=-25] (s3);
                \draw[edge,dashed] (s3) to[out=160, in=20] (s4);
                \draw[edge, dashed] (s4) to[out=-155, in=65] (srm1);
                \draw[edge] (srm1) to[out=-110, in=110] (sr);
                \draw[edge] (sr) to[out=-67.5, in=140] (vl);
                \draw[edge] (vl) to[out=-40, in=157.5] (v1);

                \draw[] (70:2.72) to[out=160, in=20] (112.5:2.72);
                \node[] (atmostk1) at (90:2.4) {\small $\geq 4k-3$};
                
                \draw[] (-157.5:2.72) to[out=-67.5, in=157.5] (-112.5:2.72);
                \node[] (atmostk1) at (-135:2.45) {\small $<k$};
                
                \draw[] (-67.5:2.72) to[out=22.5, in=-90] (0:2.72);
                \draw[] (0:2.72) to[out=90, in=-22.5] (65:2.72);
                \node[] (atleastk) at (0:2.2) {\small $\geq k$};

                \draw[dotted, edge, bend right=80] (v2) to (s1);
                \draw[dotted, edge, bend right=80] (s1) to (s2);
                \draw[dotted, edge, bend right=80] (srm1) to (sr);
                \draw[dotted, edge, bend right=80] (sr) to (v1);

              \end{tikzpicture}
          \caption{The structure of $C$ in $D$. The solid and dashed arcs represent the arcs of $C$. A dotted arc from $u$ to $v$ illustrates the existence of two internally-disjoint directed path from $u$ to $v$ in $D$.}
          \label{fig:digirth_large_C_k_k}
        \end{center}
    \end{figure}

We now define $d_i = c_{i+i_0}$ for $i=0,\dots,5$. If $i$ is an index larger than $5$, we identify $d_i$ with $d_{i \mod 6}$. For every $i=0,\dots,5$, let $X_i$ be the set of vertices reachable from $d_i$ in $D-d_{i+1}$. Similarly, if $i$ is larger than $5$, we identify $X_i$ with $X_{i \mod 6}$.

\begin{claim}\label{claim:X_i_cap_C}
    For every $i=0,\dots,5$, $X_i \cap V(C) = V(C[d_i,d_{i+1}[)$.
\end{claim}

\begin{proofclaim}
    We first consider the case $i \in \{0,\dots,4\}$. Assume for a contradiction that there exists a directed path  $P$ from $V(C[d_i,d_{i+1}[)$ to $V(C) \setminus V(C[d_i,d_{i+1}])$ with  internal vertices disjoint from $C$.
    Let $u$ be $\init(P)$ and $v$ be $\term(P)$.
    If $v \in V(C]d_{i+1},a])$, then $C[b,u] \cup P \cup C[v,a]$ is a $(b,a)$-path which avoids $d_{i+1}$, a contradiction.
    Otherwise $v \in V(C[b,d_i[)$, and then $P$ has length at least $k$ because, since $C$ is isometric,
    \[ \dist_P(u,v) \geq \dist_C(u,v) \geq \dist_C(d_5,b) \geq \dist_C(c_{i_1},a) \geq k. \]
    But then $P \cup C[u,v]$ is a subdivision of $C(k,k)$ with source $u$ and sink $v$, a contradiction to $D$ being a counterexample.

    Now suppose $i=5$. 
    Consider first a directed path $P$ from $u \in V(C[d_5,a])$ to $v \in V(C]d_0,d_5[)$ internally disjoint from $C$.
    Then $P$ has length at least $k$ because 
    \[ \dist_P(u,v) \geq \dist_C(u,v) \geq \dist_C(b,d_0) \geq k. \]
    Thus $P \cup C[u,v]$ is a subdivision of $C(k,k)$ with source $u$ and sink $v$, a contradiction.
    Consider finally a directed path $P$ from $u \in V(C[b,d_0[)$ to $v \in V(C]d_0,d_5[)$.
    Then $d_0$ is not a $(b,a)$-cut in $D$, a contradiction.
    
\end{proofclaim}

\begin{claim}\label{claim:X_i_almost_disjoint}
    For every distinct $i,j \in \{0,\dots,5\}$, we have
    \begin{enumerate}[label=$(\roman*)$]
        \item $X_i \cap X_j = \emptyset$ if $i \not\in \{j-1,j+1\}$, and
        \item $X_i \cap X_{i+1} = \emptyset$ if $v_0 \not\in X_i \cup X_{i+1}$.
    \end{enumerate}
\end{claim}

\begin{proofclaim}
    We first prove $(i)$. Let us fix two distinct integers $i,j \in \{0,\dots,5\}$ such that $i \not\in \{j-1, j+1\}$.
    Assume for a contradiction that $X_i \cap X_j \neq \emptyset$. 
    Note that $X_i \cap X_j \neq V(D)$ because $d_{i+1} \notin X_i$. Therefore, since $D$ is strongly connected by Claim~\ref{claim:strong}, there is an arc $uv$ such that $u\in X_i \cap X_j$ and $v \in V(D) \setminus (X_i \cap X_j)$.
    Assume first that $v\in X_j \setminus X_i$. Since $v\notin X_i$, we must have $v = d_{i+1}$. Hence $d_{i+1} \in X_j$, a contradiction to Claim~\ref{claim:X_i_cap_C}. 
    Symmetrically, if $v\in X_{i}\setminus X_j$ then $v = d_{j+1}$ by definition of $X_j$, implying that $d_{j+1} \in X_i$, a contradiction to Claim~\ref{claim:X_i_cap_C}. 
    Finally if $v\notin (X_i \cup X_j)$, by definition of $X_i$ and $X_j$, we must have $v = d_{i+1} = d_{j+1}$, a contradiction.
    This proves $(i)$.
    
    \medskip
    
    We now prove $(ii)$. Assume for a contradiction that $v_0\notin X_i \cup X_{i+1}$ and $X_i\cap X_{i+1} \neq \emptyset$. Recall that $X_i \cap X_{i+1} \neq V(D)$ because $d_{i+1} \notin X_i$.
    Therefore, since $D$ is strongly connected, there is an arc $uv$ such that $u \in X_i \cap X_{i+1}$ and $v\in V(D) \setminus (X_i \cap X_{i+1})$. 
    First, if $v \in V(D) \setminus (X_i \cup X_{i+1})$ then, by definition of $X_i$, $v$ must be $d_{i+1}$, 
    and by definition of $X_{i+1}$, $v$ must be $d_{i+2}$, a contradiction. 
    Next if $v \in X_i \setminus X_{i+1}$, then by definition of $X_{i+1}$, $v$ must be $d_{i+2}$, but $d_{i+2} \notin X_i$ by Claim~\ref{claim:X_i_cap_C}, a contradiction. 
    Then we may assume that $v \in X_{i+1} \setminus X_i$, and by definition of $X_i$, $v$ must be $d_{i+1}$.
    As $u \in X_{i+1}$, there is a directed path $P$ from $V(C[d_{i+1}, d_{i+2}[)$ to $u$ in $D-d_{i+2}$ internally disjoint from $C$.
    Let $x$ be $\init(P)$.
    If $x \neq d_{i+1}$, then the union of $P \cup ud_{i+1}$ (which has length at least $k$ because $\digirth(D) \geq 8k-6 \geq 2k$ and $\dist(d_{i+1},x) \leq k$ by Claim~\ref{claim:dist_ci_cj}) and $C[x, d_{i+1}]$ (which has length at least $k$ because $|V(C)| \geq 8k-6 \geq 2k$ and $\dist(d_{i+1},x) \leq k$ by Claim~\ref{claim:dist_ci_cj}) is a subdivision of $C(k,k)$ with source $x$ and sink $d_{i+1}$, a contradiction.

    So we assume that $x = d_{i+1}$, that is $P \cup ud_{i+1}$ is a cycle $C'$ with $V(C') \cap V(C) = \{d_{i+1}\}$, and $u \in V(C') \cap X_i \cap X_{i+1}$.
    Let $w$ be the vertex in $V(C') \cap X_i$ such that $\dist_{C'}(w,d_{i+1})$ is maximum (the existence of $w$ is guaranteed because $u\in V(C') \cap X_i$), and let $Q$ be a $(d_i,w)$-path in $D\ind{X_i}$.
    If $\dist_{C'}(w,d_{i+1}) \leq k-1$, then $C'[d_{i+1}, w]$ has length at least $k$,
    and $C[d_{i+1},d_i] \cup Q$  has length at least $k$. Moreover, the directed paths $C'[d_{i+1}, w]$ and  $C[d_{i+1},d_i] \cup Q$ are internally vertex-disjoint by the choice of $P$, $w$ and $Q$.
    Hence their union is a subdivision of $C(k,k)$ with source $d_{i+1}$ and sink $w$, a contradiction.
    Henceforth we suppose that $\dist_{C'}(w,d_{i+1}) \geq k$.

    We now prove the following statement.
    \begin{equation}
        \label{eq:statement_no_path_R}
        \text{$D - d_{i+1}$ does not contain any directed path $R$ from $V(C')$ to $V(C)$.}
    \end{equation}
    Assume for a contradiction that such a directed path $R$ exists. We assume that $R$ is internally disjoint from $V(C')\cup V(C)$, for otherwise we can extract a subpath of $R$ with this extra property. 
    Let $y=\init(R)$ and $z=\term(R)$.
    Then by Claim~\ref{claim:X_i_cap_C}, $z$ belongs to $V(C]d_{i+1},d_{i+2}])$.
    Observe that $y\in V(C']d_{i+1},w[)$, for otherwise $y$ belongs to $X_i$ and so does $z$, a contradiction to Claim~\ref{claim:X_i_cap_C}.
    But then the union of $R \cup C[z,d_{i+1}]$
    and $C'[y,d_{i+1}]$ is a subdivision of $C(k,k)$ with source $y$ and sink $d_{i+1}$, a contradiction to $D$ being a counterexample. This shows~(\ref{eq:statement_no_path_R}).
    
    Let $U$ be the set of vertices reachable from $d_i$ in $D\setminus d_{i+1}t$ where $t$ is the successor of $d_{i+1}$ in $C$.
    We claim that $U \subseteq X_i \cup X_{i+1}$. Let $u$ be any vertex in $U$. By definition, there is a directed path $R'$ from $d_i$ to $u$ in $D\setminus d_{i+1}t$.
    If $d_{i+1} \not\in V(R')$, then $u \in X_i$. Else if $d_{i+1} \in V(R')$ and $d_{i+2} \not\in V(R')$, then $u \in X_{i+1}$. Henceforth assume that both $d_{i+1}$ and $d_{i+2}$ belong to $R'$. Observe that $d_{i+1}$ is before $d_{i+2}$ along $R'$, otherwise $d_{i+2}\in X_i$, a contradiction to Claim~\ref{claim:X_i_cap_C}. 
    Since $d^+_D(d_{i+1}) = 2$, the successor of $d_{i+1}$ in $R'$ is also its successor in $C'$. Hence  $R'[d_{i+1},d_{i+2}]$ contains a subpath $R$ from  $V(C')\setminus \{d_{i+1}\}$ to $V(C) \setminus \{d_{i+1}\}$ internally disjoint from $V(C') \cup V(C)$, a contradiction to~(\ref{eq:statement_no_path_R}).
    
    This proves that $U \subseteq X_i \cup X_{i+1}$ and in particular, $v_0 \not\in U$.
    Set $v'_0=d_{i+1}$, $D' = D\ind{U}$. Then $D'$ equipped with $v'_0$ is such that every vertex in $U$ has out-degree $2$ in $D'$ except $v'_0$ which has out-degree at least $1$.
    By minimality of $|V(D)|$, $D'$ contains a subdivision of $C(k,k)$ and so does $D$, a contradiction.
    
\end{proofclaim}

By $(i)$ of the previous claim, there is an index $j \in \{0,\dots, 5\}$ such that $v_0 \not\in X_{j-1}$. Since $D$ is strongly connected, $\bigcup_{i=0}^5 X_i = V(D)$, and so there is an index $i \in \{0,\dots, 5\}$ such that $v_0 \in (X_i \setminus X_{i-1})$. From now on, we fix such an index $i\in \{0,\ldots,5\}$, and we set
\begin{align*}
    Y_0 &= X_{i-1} \cup X_i \cup X_{i+1} \cup X_{i+2},\\
    Y_1 &= X_{i+3}, \text{~and}\\
    Y_2 &= X_{i+4}.
\end{align*}
Note that $v_0 \not\in X_{i-1} \cup X_{i+2}$, and so,
by Claim~\ref{claim:X_i_almost_disjoint}, $Y_0,Y_1,Y_2$ are pairwise vertex-disjoint.
Moreover, $Y_0 \setminus V(C), Y_1 \setminus V(C), Y_2 \setminus V(C)$ are pairwise non adjacent by definition of the $X_j$s ({\it i.e.} there is no arc of $D$ with head and tail in different parts of  $(Y_0 \setminus V(C), Y_1 \setminus V(C), Y_2 \setminus V(C))$).
Consider a connected component $A$ of $D - V(C)$ of maximal size, that is with $|A| = \rho(C)$.
Then $A$ is included in one of $Y_0,Y_1,Y_2$. Let $j \in \{1,2\}$ be such that $A \cap Y_j = \emptyset$.
Let $q$ be the predecessor of $d_{i+j+3}$ in $C$.
Let $S$ be the set of vertices reachable from $q$ in $D-d_{i+j+3}$.
Observe that $S$ is a subset of $X_{i+j+2}$.
We claim that $D\ind{S}$ is not acyclic. Indeed, for every vertex $u\in S$, $N^+_D(u) \subseteq N^+_{D\ind{S}}(u) \cup \{d_{i+j+3}\}$. Since $v_0 \not\in S$, for every vertex $u\in S$, $d^+_{D\ind{S}}(u) \geq d^+_D(u)-1 =1$. Therefore $D\ind{S}$ has minimum out-degree at least $1$. Let $C'$ be an isometric cycle in $D\ind{S}$.

Let us show that $C'$ is an isometric cycle in $D$.
Suppose on the contrary that there is a directed path $P$ from $x \in V(C')$ to $y \in V(C')$ internally disjoint from $V(C')$ of length smaller than $\dist_{C'}(x,y)$.
As $P$ is not included in $S$, $P$ must contain $d_{i+j+3}$.
Let $d_\iota$ be the last vertex along $P$ in $\{d_\ell \mid \ell =0,\dots, 5\}$.
We have $y \in X_\iota$ by definition of $\iota$ and $y \in X_{i+j+2}$ because $y\in S \subseteq X_{i+j+2}$. Therefore $y\in X_\iota \cap X_{i+j+2}$. Since $v_0 \not\in X_{i+j+1} \cup X_{i+j+2} \cup X_{i+j+3}$, by Claim~\ref{claim:X_i_almost_disjoint},  we deduce that $\iota=i+j+2$.
Hence $P$ contains a directed path from $d_{i+j+3}$ to $d_{i+j+2}$. This implies:
\[
{
\everymath={\displaystyle}
\renewcommand{\arraystretch}{2}
\begin{array}{r l l}
    dist_{C'}(x,y) &\geq \length(P) &\text{by definition of $P$}\\
            &\geq \dist_{D}(d_{i+j+3},d_{i+j+2}) &\text{because $P$ contains $d_{i+j+3}$ and $d_{i+j+2}$} \\
            &= \dist_{C}(d_{i+j+3},d_{i+j+2}) &\text{because $C$ is isometric} \\
            &= |V(C)| - \dist_{C}(d_{i+j+2},d_{i+j+3}) &\\
            &\geq (8k-6) - (k-1) \geq k &\text{by Claim~\ref{claim:dist_ci_cj}.}\\
\end{array}
}
\]
Therefore both $P$ and $C'[x,y]$ have length at least $k$, implying that $P\cup C'[x,y]$ is a subdivision of $C(k,k)$ with source $x$ and sink $y$, a contradiction to $D$ being a counterexample.
This proves that $C'$ is isometric in $D$.

By definition, $N(A) \subseteq V(C)\cup A$. Since $D$ is strongly connected, $A$ has an in-neighbour in $V(C)$. Since $A \cap Y_j = \emptyset$ by choice of $j$, $A$ has an in-neighbour in $C[d_{i+j+3}, d_{i+j+2}[$.
Hence, the connected component in $D-V(C')$ which contains $A$ is strictly larger than $A$,
which contradicts the maximality of $\rho(C)$, and concludes the proof of the theorem.
\end{proof}

\addtocounter{theorem}{6}

\subsection{Subdivisions of \texorpdfstring{$C(2,2)$}{C(2,2)} in oriented graphs}

In this section, we improve Theorem~\ref{thm:mader_spindle_d+} when $k=2$ as follows.

\begin{theorem}\label{thm:mader_small_spindle_d+}
    Every oriented graph $D$ with $\delta^+(D) \geq 2$ contains a subdivision of $C(2,2)$.
\end{theorem}

\begin{proof}
Suppose for contradiction that there exists an oriented graph $D$ with $\delta^+(D) \geq 2$
that contains no subdivision of $C(2,2)$.
Assume that $|V(D)|$ is minimum, and that among such minimum counterexamples, $|A(D)|$ is minimum.

\begin{claim}\label{claim:out-degree-2}
    For every vertex $v \in V(D)$, $d^+(v) = 2$.
\end{claim}

\begin{proofclaim}
    If $v$ is a vertex with at least $3$ out-neighbours $w_1,w_2,w_3$, then
    $D\setminus v w_3$ is a smaller counterexample.
\end{proofclaim}

\begin{claim}\label{claim:strong}\label{claim:no-source}
    $D$ is strongly connected. In particular, $d^-(v) \geq 1$ for every vertex $v$.
\end{claim}

\begin{proofclaim}
    Let $C$ be a terminal strongly connected component of $D$. Then $C$ is also a counterexample, so by minimality of $D$ we have $D=C$, and $D$ is strongly connected.
\end{proofclaim}

\begin{claim}\label{claim:in-degree_at_least_2}
    For every vertex $v \in V(D)$, $d^-(v) \geq 2$.
\end{claim}

\begin{proofclaim}
    Suppose that $v$ is a vertex which has at most one in-neighbour. By Claim~\ref{claim:no-source}, it must have a unique in-neighbour $u$, and let $w_1,w_2$ be its two out-neighbours.
    If $w_1$ is non adjacent to $u$, then consider $D' = (D - v) \cup uw_1$.
    By minimality of $D$, $D'$  contains a subdivision $F$ of $C(2,2)$, and as $F \not\subseteq D$, we have $uw_1 \in A(F)$.
    But then $(F \setminus uw_1) \cup uv \cup vw_1 \subseteq D$ is a subdivision of $C(2,2)$.
    Hence there is an arc between $u$ and $w_1$.
    Similarly, there is an arc between $u$ and $w_2$.
    
    If $w_1u,w_2 u \in A(D)$, then the union of the directed paths $(v, w_1, u)$ and $(v,w_2,u)$ yields a copy of $C(2,2)$ in $D$.
    Moreover, if $uw_1,uw_2 \in A(D)$, then $d^+(u) \geq 3$, a contradiction to Claim~\ref{claim:out-degree-2}.
    Hence, without loss of generality, $w_1u,uw_2 \in A(D)$.
    
    As $D$ is strongly connected, there is a directed path $P$ from $w_2$ to $\{u,v,w_1\}$ with internal vertices disjoint from
    $\{u,v,w_1,w_2\}$. The terminal vertex of $P$ is not $v$, as the only in-neighbour of $v$ is $u$.
    So the terminal vertex of $P$ is either $u$ or $w_1$.
    If it is $w_1$, then the union of the directed paths $(u,v,w_1)$ and $uP$ yields a subdivision of $C(2,2)$.
    If $u$ is the end-vertex of $P$, then the union of the directed paths $(v,w_1,u)$ and $vw_2 \cup P$ yields a subdivision of $C(2,2)$.
    In both cases, we find a subdivision of $C(2,2)$ in $D$, a contradiction.
\end{proofclaim}

\begin{claim}\label{claim:2-diregular}
    $D$ is 2-diregular.
\end{claim}
\begin{proofclaim}
    By Claim~\ref{claim:out-degree-2}, we know that, for each $v\in V(D)$, $d^+(v) = 2$. It implies that $|A(D)| = \sum_{v\in V(D)} d^+(v) = 2|V(D)|$. Since $|A(D)|$ is also equal to $\sum_{v\in V(D)} d^-(v)$, we get by Claim~\ref{claim:in-degree_at_least_2} that for every vertex $v$ of $D$, $d^-(v) = d^+(v) = 2$, which implies that $D$ is 2-diregular.
\end{proofclaim}

\begin{claim}\label{claim:out-neighbour_adjacent_in-neighbour}
    For every arc $vw$, $w$ has a neighbour in $N^-(v)$.
\end{claim}

\begin{proofclaim}
    Let $u_1,u_2$ be the in-neighbours of $v$.
    If $w$ has no neighbour in $\{u_1,u_2\}$, consider $D' = D - v \cup u_1w \cup u_2 w$.
    By minimality of $D$, there exists a subdivision $F$ of $C(2,2)$ in $D'$.
    If neither $u_1w$ nor $u_2w$ belongs to $A(F)$, then $F \subseteq D$, a contradiction.
    If both $u_1w$ and $u_2w$ belong to $A(F)$, then $w$ is the sink of $F$ and $F \setminus w_1 \cup v$ is a subdivision of $C(2,2)$ in $D$,  a contradiction.
    If exactly one of $u_1w$ and $u_2w$ belongs to $A(F)$, say $u_1w$, then $F \setminus u_1w \cup \{u_1v,vw\}$ is a subdivision of $C(2,2)$ in $D$, a contradiction.
    Hence $w_1$ has a neighbour in $\{u_1,u_2\}$.
\end{proofclaim}

\begin{claim}\label{claim:neighbourhood_structure}
    For every vertex $v$ with in-neighbourhood $u_1,u_2$ and out-neighbourhood $w_1,w_2$, either $\{ w_1u_1, w_2u_2 \} \subseteq A(D)$ or
    $\{ w_1u_2, w_2u_1 \} \subseteq A(D)$. In particular, every vertex belongs to two different directed triangles.
\end{claim}

\begin{proofclaim}
    By Claim~\ref{claim:out-neighbour_adjacent_in-neighbour}, $w_1$ has a neighbour in $\{u_1,u_2\}$.
    Without loss of generality, suppose that it is $u_1$.
    We now show that $w_1u_1 \in A(D)$, so assume for a contradiction that $u_1w_1 \in A(D)$. By Claim~\ref{claim:out-neighbour_adjacent_in-neighbour} $u_1$ has an in-neighbour $x$ which is also a neighbour of $w_1$.
    
    If $x=u_2$, then $D\ind{\{u_1,u_2,v,w_1\}}$ contains a copy of $C(2,2)$ with source $u_2$ and sink $w_1$.
    If $x=w_2$, then either $w_2w_1 \in A(D)$ and $w_1$ has in-degree $3$, a contradiction to Claim~\ref{claim:2-diregular}, or
    $w_1w_2 \in A(D)$ and $D\ind{\{u_1,v,w_1,w_2\}}$ contains a copy of $C(2,2)$ with source $u_1$ and sink $w_2$.
    Hence $x,u_1,u_2,v,w_1,w_2$ are distinct.
    
    Moreover, $xw_1 \not\in A(D)$ for otherwise $w_1$ has in-degree $3$ contradicting Claim~\ref{claim:2-diregular}.
    Hence $w_1 x \in A(D)$. Consider the out-neighbour $y$ of $w_1$ distinct from $x$.
    By Claim~\ref{claim:out-neighbour_adjacent_in-neighbour}, $y$ has a neighbour in $N^-(w_1) = \{u_1,v\}$.
    If $v$ is a neighbour of $y$, then $y\in\{u_2,w_2\}$.
    If $y=w_2$, then $D\ind{\{u_1,v,w_1,w_2\}}$ contains a copy of $C(2,2)$ with source $u_1$ and sink $w_2$.
    If $y = u_2$, then the union of the directed paths $(w_1, u_2, v)$ and $(w_1, x, u_1, v)$ yields a subdivision of $C(2,2)$ with source $w_1$ and sink $v$.
    
    Hence $y$ is not a neighbour of $v$, and so $y$ is a neighbour of $u_1$.
    If $u_1y \in A(D)$, then $u_1$ has out-degree at least $3$, contradicting Claim~\ref{claim:out-degree-2}.
    If $yu_1 \in A(D)$, then $D\ind{\{w_1,x,u_1,y\}}$ contains a copy of $C(2,2)$ with source $w_1$ and sink $u_1$.
    In both cases, we reach a contradiction.
    
    This proves that $w_1u_1 \in A$.
    Similarly, $w_2$ has an out-neighbour in $N^-(v) = \{u_1,u_2\}$.
    If $w_2u_1 \in A(D)$, then $D[\{v,w_1,w_2,u_1\}]$ contains a copy of $C(2,2)$ with source
    $v$ and sink $u_1$, a contradiction.
    Hence $w_2u_1 \notin A(D)$, and so $w_2u_2 \in A(D)$ as claimed.
\end{proofclaim}

\begin{claim}\label{claim:triangle_intersection}
    Let $t_1$ and $t_2$ be two distinct directed triangles of $D$. Then $|V(t_1) \cap V(t_2)| \leq 1$.
\end{claim}
\begin{proofclaim}
    Since $D$ is an oriented graph, then it is clear that $|V(t_1) \cap V(t_2)| \leq 2$. Assume now that $V(t_1) = \{x,y,z\}$ and $V(t_2) = \{x,y,w\}$, where $z\neq w$. Assume without loss of generality that $xy\in A(t_1)$, then $xy\in A(t_2)$ because $x$ and $y$ must be adjacent in $t_2$, and $D$ does not contain any digon.
    Now $t_1\cup t_2$ contains a copy of $C(2,2)$ with source $y$ and sink $x$, a contradiction.
\end{proofclaim}

Consider the undirected auxiliary graph $H$ whose vertices are the directed triangles in $D$, and such that
two directed triangles $t$ and $t'$ of $D$ are adjacent in $H$ if and only if they share a common vertex.

By Claim~\ref{claim:triangle_intersection} and Claim~\ref{claim:2-diregular}, $H$ is a subcubic graph. Moreover, by Claim~\ref{claim:neighbourhood_structure}, $H$ must be a cubic graph.
In particular $H$ is not a forest and so it contains an induced cycle $C=(t_1, \dots, t_k,t_1)$. Recall that $t_1,\dots,t_k$ are directed triangles of $D$.
Let $t_1 = (x,y,z,x)$ and suppose (by possibly relabelling $t_1$ and $C$) that $V(t_1) \cap V(t_k) = \{x\}$ and $V(t_1) \cap V(t_2) = \{y\}$.
Let $P$ be a directed path in $D$ with vertices in $V(t_2) \cup \dots \cup V(t_k)$ from $y$ to $x$.
Observe that $z \not\in V(P)$ because $C$ is an induced cycle of $H$.
Then the union of $P$ and the path $(y,z,x)$ is a subdivision of $C(2,2)$ in $D$, a contradiction.
This proves the theorem.
\end{proof}

\subsection{Subdivisions of out-stars}
This section is devoted to the proof of Theorem~\ref{thm:out_star}.

\thmoutstar*

\begin{proof}
    Let $D$ be such a digraph. 
    By taking $|A(D)|$ minimal, we can suppose that $d^+(v)=k$ for every vertex $v\in V(D)$.
    Let $W$ be the set of vertices at distance at least $\ell$ from $u$.
    If there are $k$ vertex-disjoint $(u,W)$-paths then these directed paths have length at least $\ell$ and so they form a copy of $S_k^{+(\ell)}$.
    Otherwise, by Menger's Theorem, there is a set $S \subseteq V(D) \setminus \{u\}$ of $k-1$ vertices such that there is no $(u,W)$-path in $D-S$.
    Let $R$ be the set of vertices reachable from $u$ in $D-S$. Then every vertex in $R$ is at distance at most $\ell-1$ from $u$, so $|R| \leq \frac{k^\ell - 1}{k - 1}$.
    As $D$ has digirth at least $\frac{k^\ell - 1}{k - 1}+1$, this implies that $D\ind{R}$ is acyclic. Let $r \in R$ be a sink in $D\ind{R}$.
    Then all the out-neighbours of $r$ in $D$ are in $S$, and so $d^+(r) \leq k-1$, a contradiction.
\end{proof}

When $k=2$, we strengthen Theorem~\ref{thm:out_star} as follows.

\begin{theorem}\label{thm:small_out_star}
    Let $\ell$ be a positive integer.
    Every digraph $D$ with $\delta^+(D) \geq 2$ and $\digirth(D) \geq 2\ell$ contains a copy of $S^{+(\ell)}_2$.
\end{theorem}

\begin{proof}
    Suppose for contradiction that there is such a digraph $D$ containing no copy of $S^{+(\ell)}_2$.
    We assume $\ell \geq 2$, the result being trivial when $\ell = 1$.
    Without loss of generality, we may also assume that $d^+(v) = 2$ for every vertex $v$ in $D$.
    By considering only one terminal strongly connected component of $D$, we can also assume that $D$ is strong.
    Let $u$ be a vertex in $D$, and let $w$ be a vertex at distance exactly $\ell$ from $u$. 
    Such a vertex exists because, as $D$ is strong, $u$ has an in-neighbour, which is at distance at least $2\ell-1 \geq \ell+1$ from $u$. 

    Let us fix $P$ a shortest directed path from $u$ to $w$.
    A {\it $P$-tricot} is a sequence of pairwise vertex-disjoint directed paths $Q_1, \dots, Q_r$ (where $r$ is the size of the tricot) such that for every $i \in [r]$:
    \begin{itemize}
        \item $V(Q_i) \cap V(P) = \{\term(Q_i),\init(Q_i) \}$,
        \item $\init(Q_1)=u$,
        \item $\init(Q_{i+1}) = \pred_P(\term(Q_i))$ if $i<r$,
        \item there is an arc from $\pred_{Q_i}(\term(Q_i))$ to $V(P]\init(Q_i), \term(Q_i)[)$.
    \end{itemize}

    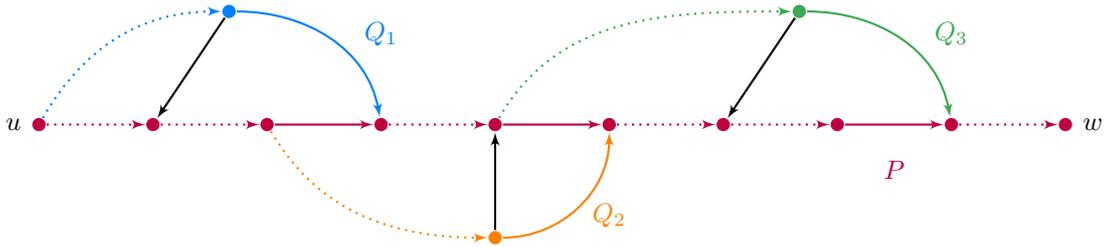
\begin{figure}[hbtp]
        \begin{center}	
              \begin{tikzpicture}[thick,scale=1, every node/.style={transform shape}]
                \tikzset{vertex/.style = {circle,fill=black,minimum size=5pt, inner sep=0pt}}
                \tikzset{littlevertex/.style = {circle,fill=black,minimum size=0pt, inner sep=0pt}}
    	    \tikzset{edge/.style = {->,> = latex'}}
                \tikzset{bigvertex/.style = {shape=circle,draw, minimum size=2em}}

                \node[vertex, purple, label=left:$u$] (u0) at (0,0){};
                \node[vertex, purple, label=right:$w$] (u9) at (13.5,0){};
                \foreach \i in {1,...,8}{
                    \node[vertex,purple] (u\i) at (\i*1.5,0) {};
                }
                \foreach \i in {0,1,3,5,6,8}{
                    \pgfmathtruncatemacro{\j}{\i+1}
                    \draw[edge, dotted, purple] (u\i) to (u\j) {};
                }
                \foreach \i in {2,4,7}{
                    \pgfmathtruncatemacro{\j}{\i+1}
                    \draw[edge, purple] (u\i) to (u\j) {};
                }
                \node[vertex, g-blue] (q1) at (2.5,1.5) {};
                \draw[edge, dotted, g-blue] (u0) to[out=60, in=180] (q1) {};
                \draw[edge, g-blue] (q1) to[out=0, in=100] (u3) {};
                \draw[edge] (q1) to (u1) {};
                \node[g-blue] at (4.5,1.2) {$Q_1$};

                \node[vertex, orange] (q2) at (6,-1.5) {};
                \draw[edge, dotted, orange] (u2) to[out=-60, in=180] (q2) {};
                \draw[edge, orange] (q2) to[out=0, in=-90] (u5) {};
                \draw[edge] (q2) to (u4) {};
                \node[orange] at (7.5,-1.2) {$Q_2$};
                
                \node[vertex, g-green] (q3) at (10,1.5) {};
                \draw[edge, dotted, g-green] (u4) to[out=60, in=180] (q3) {};
                \draw[edge, g-green] (q3) to[out=0, in=100] (u8) {};
                \draw[edge] (q3) to (u6) {};
                \node[g-green] at (12,1.2) {$Q_3$};
                
                \node[purple] at (11.25,-0.6) {$P$};
                
            \end{tikzpicture}
          \caption{An example of a $P$-tricot of size $3$. Dotted arcs represent directed paths.}
          \label{fig:P_tricot_1}
        \end{center}
    \end{figure}

    Let us first prove that $D$ admits a $P$-tricot. Let $Q$ be a maximum directed path in $D$, starting on $u$, that is disjoint from $V(P)\setminus \{u\}$. Since $Q$ is maximum, the two out-neighbours of $\term(Q)$ belong to $V(P)\cup V(Q)$. We have that the length of $Q$ is at most $\ell-1$, for otherwise the union of $P$ and $Q$ contains a copy of $S_2^{+(\ell)}$, a contradiction to the choice of $D$.  This implies that the out-neighbourhood of $\term(Q)$ is in $V(P)\setminus V(Q)$, for otherwise $D\ind{V(Q)}$ contains a directed cycle of length at most $\ell-1$, a contradiction. Let $x$ be the out-neighbour of $\term(Q)$ which is the furthest from $u$ and $y$ be its other out-neighbour. Let $Q'$ be the extension of $Q$ with $x$, then $(Q')$ is a $P$-tricot of size one since $Q'$ starts at $u$, intersects $P$ exactly on $\{u,x\}$ and $y$ belongs to $P[u,x[$.
    
    Among all $P$-tricots of $D$, we choose one with maximum size $r$ and denote it by $\mathcal{T}$.
    Let $P_1$ be the directed path corresponding to the concatenation of $Q_1$ and all $P[\term(Q_{i-2}),\init(Q_i)]\cdot Q_i$ for odd $i\in \{3,\ldots,r\}$. 
    Let $P_2$ be the directed path corresponding to the concatenation of all $P[\term(Q_{i-2}),\init(Q_i)]\cdot Q_i$ for even $i\in \{2,\ldots,r\}$ (we identify $\term(Q_0)$ with $u$).

    Observe that $P_1$ and $P_2$ are two directed paths starting from $u$, and that they intersect exactly on $\{u\}$. Note that $P_1$ can be completed into a $(u,w)$-path $\Tilde{P}_1 = P_1 \cdot P[\term(P_1),w]$ disjoint from $P_2-\term(P_2)$. Since it is a $(u,w)$-path, $\Tilde{P}_1$ has length at least $\ell$. Therefore $P_2$ has length at most $\ell$, for otherwise $\Tilde{P}_1 \cup (P_2-\term(P_2))$ is a copy of $S^{+(\ell)}_2$ in $D$.
    Analogously, $P_2$ can be completed into a $(u,w)$-path, which implies that $P_1$ has length at most $\ell$.
    
    Let $i\in \{1,2\}$ be such that $P_i$ does not contain $Q_r$ and let $j\in \{1,2\}$ be different from $i$. Let $v$ be $\pred_P(\term(Q_r))$. We consider $Q'$ a maximal directed path starting from $v$ in $D - (V(P) \cup V(P_1) \cup V(P_2) \setminus \{ v \})$. Let $t$ be $\term(Q')$, let $P_i'$ be the concatenation $P_i \cdot P[\term(P_i),v] \cdot Q'$ and $P_j'$ be the concatenation $P_j \cdot P[\term(Q_r),w]$. See Figure~\ref{fig:tricot_Pi_Pj} for an illustration.  
        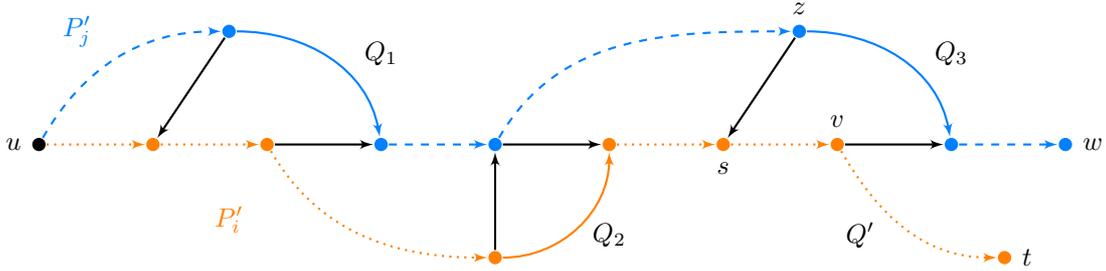
\begin{figure}[hbtp]
        \begin{center}	
              \begin{tikzpicture}[thick,scale=1, every node/.style={transform shape}]
                \tikzset{vertex/.style = {circle,fill=black,minimum size=5pt, inner sep=0pt}}
                \tikzset{littlevertex/.style = {circle,fill=black,minimum size=0pt, inner sep=0pt}}
    	    \tikzset{edge/.style = {->,> = latex'}}
                \tikzset{bigvertex/.style = {shape=circle,draw, minimum size=2em}}

                \node[vertex, label=left:$u$] (u0) at (0,0){};
                \node[vertex, g-blue, label=right:$w$] (u9) at (13.5,0){};
                \foreach \i in {1,2,5}{
                    \node[vertex, orange] (u\i) at (\i*1.5,0) {};
                }
                \foreach \i in {3,4,8}{
                    \node[vertex, g-blue] (u\i) at (\i*1.5,0) {};
                }
                \node[vertex, orange, label=above:$v$] (u7) at (10.5,0) {};
                \node[vertex, orange, label=below:$s$] (u6) at (9,0) {};
                
                \foreach \i in {3,8}{
                    \pgfmathtruncatemacro{\j}{\i+1}
                    \draw[edge, dashed, g-blue] (u\i) to (u\j) {};
                }
                \foreach \i in {0,1,5,6}{
                    \pgfmathtruncatemacro{\j}{\i+1}
                    \draw[edge, dotted, orange] (u\i) to (u\j) {};
                }
                \foreach \i in {2,4,7}{
                    \pgfmathtruncatemacro{\j}{\i+1}
                    \draw[edge] (u\i) to (u\j) {};
                }
                \node[vertex, dashed, g-blue] (q1) at (2.5,1.5) {};
                \draw[edge, dashed, g-blue] (u0) to[out=60, in=180] (q1) {};
                \draw[edge, g-blue] (q1) to[out=0, in=100] (u3) {};
                \draw[edge] (q1) to (u1) {};

                \node[] at (4.5,1.2) {$Q_1$};
                \node[] at (7.5,-1.2) {$Q_2$};
                \node[] at (10.8,-1.2) {$Q'$};

                \node[vertex, orange] (q2) at (6,-1.5) {};
                \node[vertex, orange, label=right:$t$] (t) at (12.7,-1.5) {};
                \draw[edge, dotted, orange] (u2) to[out=-60, in=180] (q2) {};
                \draw[edge, orange] (q2) to[out=0, in=-90] (u5) {};
                \draw[edge] (q2) to (u4) {};
                \draw[edge, dotted, orange] (u7) to[out=-60, in=180] (t) {};
                
                \node[orange] at (2.5,-1) {$P_i'$};
                \node[g-blue] at (0.5,1.5) {$P_j'$};
                
                \node[vertex, g-blue, label=above:$z$] (q3) at (10,1.5) {};
                \draw[edge, dashed, g-blue] (u4) to[out=60, in=180] (q3) {};
                \draw[edge, g-blue] (q3) to[out=0, in=100] (u8) {};
                \draw[edge] (q3) to (u6) {};
                \node[] at (12,1.2) {$Q_3$};
            \end{tikzpicture}
          \caption{An illustration of the paths $P_i'$ and $P_j'$.}
          \label{fig:tricot_Pi_Pj}
        \end{center}
    \end{figure}
    
    Since $P'_j$ is a $(u,w)$-path, $P_j'$ has length at least $\ell$. If also $P_i'$ has length at least $\ell$, then $P_i'$ and $P_j'$ are two directed paths of length at least $\ell$, sharing their source $u$ and no other vertices, a contradiction since $D$ is a counterexample. We know that $t$ has two out-neighbours $x,y$ in $D$. Since $Q'$ is a maximal directed path, we know that both $x$ and $y$ belong to $V(P_i') \cup V(P_j') \cup V(P)$. 
    \begin{itemize}
        \item First if one of $x,y$, say $x$, belongs to $V(P_i')$, $P_i'[x,t]  \cup \term(P_i')x$ is a directed cycle of length at most
         \[ 
         \length(P_i'[x,t]) + 1 \leq \length(P_i') + 1 \leq \ell, 
         \] 
         where in the last inequality we used that $P_i'$ has length at most $\ell - 1$. This is a contradiction since $\digirth(D) \geq 2\ell > \ell - 1$.
        
        \item Else if one of $x,y$, say $x$, belongs to $V(P[u,v])$, then  $P[x,v] \cdot Q' \cup tx$ is a directed cycle of length at most:
        \[
        {
        \everymath={\displaystyle}
        \renewcommand{\arraystretch}{2}
        \begin{array}{r l l}
            &\length(P[x,v]) + \length(Q'[v,t]) + 1 &\\
            &\leq \length(P[u,v]) + \length(Q'[v,t]) + 1 &\text{~~because $u$ is before $x$ in $P$}\\
            &\leq \length(P_i'[u,v]) + \length(Q'[v,t]) + 1 &\text{~~because $P$ is a shortest path}\\
            &=\length(P_i') + 1 \leq \ell,
        \end{array}
        }
        \]
        a contradiction since $\digirth(D) \geq 2\ell > \ell$.
        
        \item Else if one of $x,y$, say $x$, belongs to $V(P_j) \setminus V(P[u,v])$, let $z$ be $\pred_{P_j}(\term(P_j))$, which is also $\pred_{Q_r}(\term(Q_r))$. By definition of $\mathcal{T}$, $z$ has an out-neighbour $s$ in $V(P]\init(Q_r), \term(Q_r)[)$. Then $P_j[x,z] \cdot zs \cdot P[s,v] \cdot Q' \cup tx$ is a directed cycle with length at most
        \[
        {
        \everymath={\displaystyle}
        \renewcommand{\arraystretch}{2}
        \begin{array}{r l l}
            &\length(P_j[x,z]) + \length(P[s,v]) + \length(Q') + 2 &\\
            &\leq \length(P_j) + \length(P[s,v]) + \length(Q') &\text{~~because $x\neq u$ and $z\neq \term(P_j)$}\\  
            &\leq \length(P_j) + \length(P_i'[u,v]) + \length(Q') &\text{~~because $P$ is a shortest path}\\
            &\leq \length(P_j) + \length(P_i') \leq 2\ell - 1, &
        \end{array}
        }
        \]
        a contradiction since $\digirth(D) \geq 2 \ell$.
        
        \item Finally if both $x$ and $y$ belong to $P[\term(Q_r),w]$, we can assume that $x$ is before $y$ on the path $P$. But then the $P$-tricot $(Q_1,\dots,Q_r,Q'\cdot (t,y))$ contradicts the maximality of $\mathcal{T}$.
    \end{itemize}
\end{proof}

\section{Further research directions}
\label{sec:conclusion}

In this work, for a fixed digraph $F$, we give some sufficient conditions on a digraph $D$ to ensure that $D$ contains $F$ as a subdivision. This is just the tip of the iceberg and many open questions arise.
In particular, the exact value of $\mader_{\dic}(F)$ is known only for very few digraphs $F$. The smallest digraph $F$ for which it is unknown is $\bid{K_3}$.

\begin{conjecture}[Gishboliner, Steiner, and Szabó~\cite{gishboliner_dichromatic_2022}]
\[
\mader_{\dic}(\bid{K_3}) = 4
\]
\end{conjecture}

In the first part of this paper, we looked for paths and cycles in large dicritical digraphs. In particular, we proved in Theorem~\ref{thm:dicritical_nolong_directed_path} that for every integer $k \geq 3$, there are infinitely many $k$-dicritical digraphs without any directed path on $3k+1$ vertices. Conversely Bermond et al.~\cite{BGHS81} proved that every connected digraph with $\delta^+(D)\geq k$ and $\delta^-(D)\geq \ell$ contains a directed path of order at least $\min\{n,k+\ell+1\}$. As every vertex in a $k$-dicritical digraph has in- and out-degree at least $k-1$, we obtain that there are finitely many $k$-dicritical digraphs with no directed path on $2k-1$ vertices. The following problem then naturally comes.

\begin{problem}
    For every integer $k\geq 3$, find the largest integer $f(k) \in [2k-1,3k]$ such that the set of $k$-dicritical $\ori{P_{f(k)}}$-free digraphs is finite.
\end{problem}

Given a digraph $F$, we say that $F$ is {\it $\delta^+$-maderian} if there is an integer $k$ such that every digraph $D$ with $\delta^+(D) \geq k$
contains a subdivision of $F$. The smallest such integer $k$ is then denoted by $\mader_{\delta^+}(F)$.
The problem of characterising $\delta^+$-maderian digraphs is widely open.
In particular, Mader~\cite{maderCOMB5} conjectured that every acyclic digraph is $\delta^+$-maderian, but this remains unproven albeit many effort to prove or disprove it (see~\cite{lochetJCTB134} for a partial answering to the conjecture).

\medskip

In the remaining of the paper, we focus on digraphs of large digirth. Given a digraph $F$, and an integer $g$, we can define $\mader_{\delta^+}^{(g)}(F)$ to be the smallest integer $k$, if it exists, such that every digraph $D$ with $\delta^+(D) \geq k$ and $\digirth(D) \geq g$ contains a subdivision of $F$. 

It is interesting to note that there are digraphs which are not $\delta^+$-maderian even when restricted to digraphs of large digirth.
Indeed, for every integers $g,d$ there is a digraph $D$ with $\digirth(D) \geq g$ and $\delta^+(D) \geq d$ such that $D$ does not contains any subdivision of $\bid{K_3}$.
Such a digraph $D$ can be easily obtained from a construction by DeVos et al.~\cite{devos2012immersing} of digraphs with arbitrarily large out-degree in which every directed cycle has odd length, by removing a few arcs in order to increase the digirth. Since every subdivision of $\bid{K_3}$ has an even directed cycle, such a digraph does not contain $\bid{K_3}$ as a subdivision.

In Theorem~\ref{thm:mader_spindle_d+}, we proved that $\mader_{\delta^+}^{(8k-4)}(C(k,k)) \leq 2$.
On the other hand, the value $8k-4$ can not be replaced by $k-1$. To see this, consider the digraph with vertex set $\mathbb{Z}/(2k-1)\mathbb{Z}$ and arc set $\{(i,i+1), (i,i+2) \mid i \in V(G)\}$.
Since $|V(D)| < 2k$, $D$ has no subdivision of $C(k,k)$.
Since $\delta^+(D)=2$ and $\digirth(G) = k-1$, we deduce that $\mader_{\delta^+}^{(k-1)}(C(k,k))>2$.
Thus the following problem arises.
\begin{problem}
    Find the minimum $g \in [k,8k-4]$ such that $\mader_{\delta^+}^{(g)}(C(k,k)) \leq 2$.
\end{problem}

In this paper, we studied the value of $\min\{\mader^{(g)}_{\dic}(X) \mid g \geq 0\}$ given a digraph $X$.
Actually, we believe that this value is upper bounded by a function of the maximum degree.

\begin{conjecture}\label{conj:mader_infty_Delta}
    There is a function $f$ such that
    for every digraph $F$ with maximum degree $\Delta$,
    there is an integer $g$ such that $\mader_{\dic}^{(g)}(F) \leq f(\Delta)$.
\end{conjecture}

This is motivated by the following result by Mader~\cite{mader1998topological}, which is somehow the undirected analog of Conjecture~\ref{conj:mader_infty_Delta}.
\begin{theorem}[Mader~\cite{mader1998topological}]
    There is a function $f$ such that
    for every graph $F$, for every graph $G$ with $\delta(G) \geq \max\{\Delta(F),3\}$, if
    $\girth(G) \geq f(F)$
    then $G$ contains a subdivision of $F$.
\end{theorem}
It was later proved by K{\"u}hn and Osthus~\cite{kuhn2002topological} that one can take $f(H) = 166 \, \frac{\log|V(H)|}{\log \Delta(H)}$, which is optimal up to the constant factor.

Harutyunyan and Mohar~\cite{HARUTYUNYAN20121823} proved that there is a positive constant $c$ such that for every large enough $\Delta,g$, there is a digraph $D$
with $\girth(D) \geq g$, $\Delta(D) \leq \Delta$ and $\dic (D) \geq c \cdot \frac{\Delta}{\log \Delta}$.
This is a generalisation to the directed case of a classical result by Bollob{\'a}s~\cite{bollobas1978chromatic}.
This implies that any function $f$ satisfying Conjecture~\ref{conj:mader_infty_Delta} is such that $f(\Delta) \geq c \cdot \frac{\Delta}{\log \Delta}$.
We are inclined to believe that this is optimal and that there is such a function $f$~in~$O\left(\frac{\Delta}{\log \Delta}\right)$.

\section*{Acknowledgement}
We are grateful to Frédéric Havet for stimulating discussions and valuable advice.

\bibliographystyle{abbrv}
\bibliography{biblio}

\end{document}